\documentclass[a4paper,12pt]{article}
\usepackage{amssymb}

\newcommand{\R}{\mathbb{R}}

\newcommand{\N}{\mathbb{N}}

\newcommand{\cuad}{{\sqcap\kern-.68em\sqcup}}

\newcommand{\norm}[1]{\|#1\|}
\newcommand{\equ}[1]{(\ref{#1})}

\newtheorem{definition}{Definition}[section]
\newtheorem{theorem}{Theorem}[section]
\newtheorem{proposition}{Proposition}[section]

\newtheorem{lemma}{Lemma}[section]
\newtheorem{corollary}{Corollary}[section]
\newtheorem{remark}{Remark}[section]
\newcommand{\bremark}{\begin{remark} \em}
\newcommand{\eremark}{\end{remark} }

\begin{document}

\begin{center}{\bf  \Large   Boundary blow-up solutions to  fractional elliptic  \medskip

equations in a measure framework}\medskip
%%%%%%%%%%%%%%%%%%%%%%%%%%%%%%%%%%%%%%%%%%%%%%%%%%%%%%%%%%%%%%%%%%%%%%
%%%%%%%%%%%%%%%%%%%%%%%%%%%%%%%%%%%%%%%%%%%%%%%%%%%%%%%%%%%%%%%%%%%%%%

 \bigskip

\medskip
{\bf Huyuan Chen\footnote{hc64@nyu.edu}\quad\ Hichem Hajaiej\footnote{hh62@nyu.edu} \quad\ Ying Wang\footnote{ywang@dim.uhicle.cl} }

\bigskip

\begin{abstract}
Let $\alpha\in(0,1)$, $\Omega$ be a bounded open domain in $\R^N$ ($N\ge 2$) with $C^2$ boundary $\partial\Omega$ and $\omega$ be the Hausdorff measure on $\partial\Omega$.
  We denote by $\frac{\partial^\alpha \omega}{\partial \vec{n}^\alpha}$  a measure
$$\langle\frac{\partial^\alpha \omega}{\partial \vec{n}^\alpha},f\rangle=\int_{\partial\Omega}\frac{\partial^\alpha f(x)}{\partial \vec{n}_x^\alpha}
d\omega(x),\quad f\in C^1(\bar\Omega),$$
where $\vec{n}_x$ is the unit outward normal vector  at point $x\in\partial\Omega$.
In this paper, we prove that  problem
 \begin{equation}\label{0.1}
 \arraycolsep=1pt
\begin{array}{lll}
 (-\Delta)^\alpha  u+g(u)=k\frac{\partial^\alpha \omega}{\partial \vec{n}^\alpha}\quad & {\rm in}\quad \bar\Omega,\\[2mm]
 \phantom{   (-\Delta)^\alpha  +g(u)}
u=0\quad & {\rm in}\quad \Omega^c
\end{array}
\end{equation}
admits a unique weak solution $u_k$ under the hypotheses that
 $k>0$, $(-\Delta)^\alpha$ denotes the fractional Laplacian with $\alpha\in(0,1)$ and $g$ is a nondecreasing function satisfying
 extra conditions. We prove that the weak solution of (\ref{0.1}) is a classical solution of
$$
 \arraycolsep=1pt
\begin{array}{lll}
\ \ \  (-\Delta)^\alpha   u+g(u)=0\quad & {\rm in}\quad  \Omega,\\[2mm]
\phantom{------\ }
\ \ \ u=0\quad & {\rm in}\quad  \R^N\setminus\bar\Omega,\\[2mm]
 \phantom{}
\lim_{x\in\Omega,x\to\partial\Omega}u(x)=+\infty.
\end{array}
$$

\end{abstract}

\end{center}
%\tableofcontents \vspace{1mm}
  \vspace{1mm}


\vspace{2mm}
%%%%%%%%%%%%%%%%%%%%%%%%%%%%%%%%%%%%%%%%%%%%%%%%%%%%%%%%%%%%%%%%%%%%%%%%%%%%%%%%%%%%%%%%%%%%%%%%%%%%%%%%%%%%%%%%%%%%%%%%%%
%%%%%%%%%%%%%%%%%%%%%%%%%%%%%%%%%%%%%%%%%%%%%%%%%%%%%%%%%%%%%%%%%%%%%%%%%%%%%%%%%%%%%%%%%%%%%%%%%%%%%%%%%%%%%%%%%%%%%%%%%%

\setcounter{equation}{0}
\section{Introduction}

\subsection{motivation }
Let $\Omega$ be a bounded open domain in $\R^N$ ($N\ge 2$) with $C^2$ boundary $\partial\Omega$. The pioneering works \cite{K,O}
obtained that  the nonlinear reaction diffusion equation
\begin{equation}\label{1.1.2}
 \arraycolsep=1pt
\begin{array}{lll}
 -\Delta u+h(u)=0\ \ \ \ &
 \mbox{in}\ \ \ &\Omega,\\[2mm]
 \phantom{ -\Delta +h(u)}
u=+\infty\ \quad&  \mbox{on}\ \ \ &\partial\Omega
\end{array}
\end{equation}
admits a solution if $h$ is a locally Lipschitz continuous function which is increasing and satisfies   Keller-Osserman condition
$$
\int_1^{+\infty}\left[\int_0^{s} h(t)dt\right]^{-\frac12} ds <+\infty.
$$
Great interests in existence, uniqueness and asymptotic behavior of boundary blow-up solution to (\ref{1.1.2}) have been taken, see \cite{BM,DG,DZZ,GLL,MV0,MV,GZ}.
 It is well known that when $h(s)=s^p$ with $p>1$, (\ref{1.1.2}) has a unique solution with boundary asymptotic behavior $\rho^{-\frac2{p-1}}(x)$,
 where $\rho(x)={\rm dist}(x,\partial\Omega)$.

Comparing with the Laplacian case, a much richer structure for the solutions set appears for the non-local case.
Recently, the authors in \cite{CFQ} obtained very different phenomena of the boundary blow-up solutions to elliptic equations
  involving the fractional Laplacian, precisely,
\begin{eqnarray}\label{1.3}
(-\Delta)^\alpha u+|u|^{p-1}u&=&0\ \quad
\hbox{in} \quad \Omega,
\nonumber\\
\label{ecbir2a}
u&=&0\quad \ \hbox{in} \quad  \Omega^c,\\
\nonumber
\label{ecbir3a}
\lim_{x\in\Omega, x\to\partial\Omega}u(x)&=&+\infty,
\end{eqnarray}
where $p>0$ and the fractional Laplacian
$(-\Delta)^\alpha $ with $\alpha\in(0,1)$  is defined by
$$(-\Delta)^\alpha  u(x)=\lim_{\epsilon\to0^+} (-\Delta)_\epsilon^\alpha u(x),$$
here for $\epsilon>0$,
$$
(-\Delta)_\epsilon^\alpha  u(x)=-\int_{\R^N\setminus B_\epsilon(x)}\frac{ u(z)-
u(x)}{|z-x|^{N+2\alpha}}  dz
$$

 The existence of boundary blow-up solution of (\ref{1.3}) is derived by  constructing appropriate super and sub-solutions and
 this construction involves the one dimensional truncated Laplacian of power functions given by
\begin{equation}\label{3.1.4}
C(\tau)=\int^{+\infty}_{0}\frac{\chi_{(0,1)}(t)|1-t|^{\tau}+(1+t)^{\tau}-2}{t^{1+2\alpha}}dt,
\end{equation}
where $\tau\in (-1,0)$ and  $\chi_{(0,1)}$ is the characteristic function of the interval $(0,1)$.
It is known that there exists a unique zero point of (\ref{3.1.4}) in $ (-1,0)$, denoting $\tau_0(\alpha)$.
Then
\begin{proposition}\cite[Theorem 1.1]{CFQ}\label{pr 0}
\label{th 3.1.1}
Assume that $\Omega$ is an open, bounded and connected domain of
class $C^2$ and $\alpha\in(0,1)$. Then we have:

\noindent
{\bf Existence}: Assume that
$$
1+2\alpha<p<1-\frac{2\alpha}{\tau_0(\alpha)},
$$
the equation \equ{1.3} possesses at least one solution
 $u$ satisfying
\begin{equation}\label{3.1.5}
0<\liminf_{x\in\Omega,x\to\partial\Omega}
u(x)d(x)^{\frac{2\alpha}{p-1}}\le\limsup_{x\in\Omega,x\to\partial\Omega}
u(x)d(x)^{\frac{2\alpha}{p-1}} <+\infty.
\end{equation}
{\bf Uniqueness}:     $u$ is the unique
solution of \equ{1.3} satisfying \equ{3.1.5}.

\noindent
{\bf Nonexistence}:  In the following three cases:
\begin{itemize}
\item[i)] For any $\tau\in(-1,0)\setminus\{-\frac{2\alpha}{p-1},\
\tau_0(\alpha)\}$ and
$$1+2\alpha<p<1-\frac{2\alpha}{\tau_0(\alpha)}\quad \mbox{ or}$$
\item[ii)] For any $\tau\in(-1,0)$ and
$$
p\ge 1-\frac{2\alpha}{\tau_0(\alpha)}\ \mbox{ or}
$$
\item[iii)] For any $\tau\in(-1,0)\setminus\{\tau_0(\alpha)\}$ and
$$
1<p\le 1+2\alpha,
$$
\end{itemize}
problem  \equ{1.3} does not have a solution $u$ satisfying
\begin{equation}\label{1.7}
0<\liminf_{x\in\Omega,x\to\partial\Omega}
u(x)d(x)^{-\tau}\le\limsup_{x\in\Omega,x\to\partial\Omega}
u(x)d(x)^{-\tau} <+\infty.
\end{equation}

\noindent
{\bf Special existence for $\tau=\tau_0(\alpha)$.}
 Assume
that
$$
\max\{1-\frac{2\alpha}{\tau_0(\alpha)}+\frac{\tau_0(\alpha)+1}{\tau_0(\alpha)},1\}<p<1-\frac{2\alpha}{\tau_0(\alpha)}.
$$
Then for any $t>0$, there is  a  positive solution $u$ of equation
(\ref{1.3})  satisfying
$$ \lim_{x\in\Omega,x\to\partial\Omega}u(x)d(x)^{-\tau_0(\alpha)}=t.$$
\end{proposition}

There are some challenging questions to ask:\\[1mm]
\emph{1. Could $\tau_0(\alpha)$ be expressed explicitly?\\[1mm]
2. With what condition of general nonlinearity makes existence hold?\\[1mm]
3. The uniqueness and nonexistence restricts in the class functions (\ref{3.1.5}) and (\ref{1.7}), so are there some solutions breaking
the assumption (\ref{3.1.5})?\\[1mm]
%4. For the special existence of boundary blow solutions, could it be extended bigger region of $p$?
}

Our interest in this article is  to introduce a new method to study the boundary blow-up solutions of semilinear fractional
elliptic equations and answer above questions. The main idea is to find suitable type measure concentrated on the whole boundary
and then by making basic estimates to prove that the corresponding weak solution solves \equ{1.3}.  Our first result is stated as follows:
\begin{proposition}\label{pr 00}
Let $\alpha\in(0,1)$ and $\tau_0(\alpha)$ is the zero point of $C(\cdot)$ when $C(\cdot)$ given by (\ref{3.1.4}), then
$$\tau_0(\alpha)=\alpha-1.$$
\end{proposition}

We observe that the critical value $1-\frac{2\alpha}{\tau_0(\alpha)}$ in Proposition \ref{pr 0} turns out to be
$\frac{1+\alpha}{1-\alpha}$.  In what follows, we would like to show the details of our new method and answer the second and third questions
in the following.

\subsection{A new method and main results}

Let $\alpha\in(0,1)$   and $\omega$ be the Hausdorff measure on $\partial\Omega$.
  We denote by $\frac{\partial^\alpha \omega}{\partial \vec{n}^\alpha}$  a measure
$$\langle\frac{\partial^\alpha \omega}{\partial \vec{n}^\alpha},f\rangle=\int_{\partial\Omega}\frac{\partial^\alpha f(x)}{\partial \vec{n}_x^\alpha}d\omega(x),\qquad f\in C^\alpha(\bar\Omega),$$
where $\vec{n}_x$ is the unit inward normal vector of $\partial\Omega$ at point $x$ and
$$\frac{\partial^\alpha f(x)}{\partial \vec{n}_x^\alpha}=\lim_{t\to0^+}\frac{f(x+t\vec{n}_x)-f(x)}{t^\alpha}.$$
In this paper,  we are concerned with the
existence and uniqueness of weak solution  to the semilinear fractional elliptic
problem
\begin{equation}\label{eq 1.1}
 \arraycolsep=1pt
\begin{array}{lll}
 (-\Delta)^\alpha  u+g(u)=k\frac{\partial^\alpha \omega}{\partial \vec{n}^\alpha}\quad & {\rm in}\quad\bar\Omega,\\[2mm]
 \phantom{   (-\Delta)^\alpha  +g(u)}
u=0\quad & {\rm in}\quad \bar\Omega^c,
\end{array}
\end{equation}
where $k>0$ and $g:\R_+\to\R_+$ is continuous.

In \cite{CH}, the authors studied problem (\ref{eq 1.1}) replaced $\frac{\partial^\alpha \omega}{\partial \vec{n}^\alpha}$ by
$\frac{\partial^\alpha \nu}{\partial \vec{n}^\alpha}$ where $\nu$ is a Radon measure concentrated on boundary measure. They proved that
such a problem has a unique weak solution if $g$ is a continuous nondecreasing function satisfying
$g(0)\ge0$ and
\begin{equation}\label{CH.g1}
 \int_1^\infty g(s)s^{-1-\frac{N+\alpha}{N-\alpha}}ds<+\infty.
\end{equation}
Moreover, \cite{CH} analyzed the isolated singularity of weak solution of (\ref{eq 1.1}) in the case that $\nu=\delta_{x_0}$ with $x_0\in\partial\Omega$.
Our aim in this article is to investigate how the Hausdorff measure on $\partial\Omega$ works on the weak solution of (\ref{eq 1.1}).

Before starting our main theorems we make precise the notion of weak solution used in this note.
\begin{definition}\label{weak solution GV}
We say that $u$ is a weak solution of (\ref{eq 1.1}), if $u\in
L^1(\Omega)$, $g(u)\in L^1(\Omega,\rho^\alpha dx)$  and
$$
\int_\Omega [u(-\Delta)^\alpha\xi+ g(u)\xi]dx=k\int_{\partial\Omega}\frac{\partial^\alpha\xi(x)}{\partial\vec{n}_x^\alpha}d\omega(x),\qquad \forall\xi\in \mathbb{X}_\alpha.
$$
where $\rho(x)={\rm dist}(x,\partial\Omega)$ and
$\mathbb{X}_{\alpha}\subset C(\R^N)$ denotes the space of functions
$\xi$ satisfying:\smallskip
\begin{itemize}
\item[]
\begin{enumerate}\item[$(i)$]
${\rm supp}(\xi)\subset\bar\Omega$;
\end{enumerate}
\begin{enumerate}\item[$(ii)$]
 $(-\Delta)^\alpha\xi(x)$ exists for all $x\in \Omega$
and $|(-\Delta)^\alpha\xi(x)|\leq C$ for some $C>0$;
\end{enumerate}
\begin{enumerate}\item[$(iii)$]
there exist $\varphi\in L^1(\Omega,\rho^\alpha  dx)$
and $\varepsilon_0>0$ such that $|(-\Delta)_\varepsilon^\alpha\xi|\le
\varphi$ a.e. in $\Omega$ for all
$\varepsilon\in(0,\varepsilon_0]$.
\end{enumerate}
\end{itemize}
\end{definition}

   Now we are ready to state our first result for problem (\ref{eq 1.1}).

\begin{theorem}\label{teo 1}
Assume that   $k>0$, $\rho(x)={\rm dist}(x,\partial\Omega)$   and $g$ is a continuous nondecreasing function satisfying
$g(0)\ge0$ and
\begin{equation}\label{g1}
 \int_1^\infty g(s)s^{-1-\frac{1+\alpha}{1-\alpha}}ds<+\infty.
\end{equation}
Then
 $(i)$ problem (\ref{eq 1.1}) admits a unique positive weak solution $u_k$;

$(ii)$ the mapping $k\to u_k$ is increasing and there exists $c_1\ge1$ independent of $k$ such that
\begin{equation}\label{1.33}
  \frac{k}{c_1}\rho(x)^{\alpha-1} \le u_k(x)\le c_1k\rho(x)^{\alpha-1},\qquad \forall x\in\Omega;
\end{equation}

$(iii)$ if we assume additionally that $g$ is $C^{\beta}$ locally in $\R$ with $\beta>0$, then
$u_k$ is a classical solution of
\begin{equation}\label{1.1}
 \arraycolsep=1pt
\begin{array}{lll}
\ \  (-\Delta)^\alpha   u+g(u)=0\quad & {\rm in}\quad  \Omega,\\[2mm]
\phantom{------\ }
\ \ u=0\quad & {\rm in}\quad  \R^N\setminus\bar\Omega,\\[2mm]
 \phantom{}
\lim_{x\in\Omega,x\to\partial\Omega}u(x)=+\infty.
\end{array}
\end{equation}

\end{theorem}

We remark that in Theorem \ref{teo 1} extends the special existence of boundary blow up solutions to
fractional elliptic equation (\ref{1.1}) with general nonlinearity $g$ in integral subcritical case with the
critical exponents $\frac{1+\alpha}{1-\alpha}$, which is larger than $\frac{N+\alpha}{N-\alpha}$.
Specially, letting $g\equiv0$, there exists infinitely many boundary blow up $\alpha-$harmonic functions.

Since $\alpha-1>-\frac{2\alpha}{p-1}$,  so  we may call  the solutions of (\ref{1.1}) as the weak boundary blow-up solution
from the asymptotic behavior (\ref{1.33}). Our second interest is to consider  the limit of weak boundary blow-up
solutions.

\begin{theorem}\label{teo 2}
  Let $g(s)=s^p$ with $p\in(0,\frac{1+\alpha}{1-\alpha})$ and $u_k$ be the weak solution of
(\ref{eq 1.1}), then

$(i)$ if $p\in(1+2\alpha ,\frac{1+\alpha}{1-\alpha})$,
then the limit of $\{u_k\}$ as $k\to\infty$ exists, denoting $u_\infty$, which is a classical solution of (\ref{1.3}).
Moreover, $u_\infty$ satisfies
\begin{equation}\label{b k 01}
\frac1{c_2}\rho(x)^{-\frac{2\alpha}{p-1}} \le u_\infty(x)\le c_2\rho(x)^{-\frac{2\alpha}{p-1}},\qquad \forall x\in\Omega,
\end{equation}
where  $c_2\ge1$.

$(ii)$ if  $p\in(0,1+ 2\alpha]$, then
$$\lim_{k\to\infty}u_k(x)=+\infty,\qquad \forall x\in \Omega.$$
\end{theorem}

We notice that  the limit of weak boundary blow-up solutions is the solution of (\ref{1.1.2})
with behavior (\ref{3.1.5}) when $p\in(1+2\alpha ,\frac{1+\alpha}{1-\alpha})$ stated in Proposition \ref{pr 0}.
As a consequence of Theorem \ref{teo 2} $(ii)$,   $p\in(0,1+ 2\alpha]$, there is no solution $u$ of (\ref{1.1.2})
such that
$$\lim_{x\in\Omega,x\in\Omega} u(x)\rho^{1-\alpha}(x)=0\ \  {\rm or}\ \ +\infty.$$

From \cite{CH} and Theorem \ref{teo 1}, the Dirac mass and Hausdorff measure have different contribution
to the solution of
$$ (-\Delta)^\alpha  u+g(u)=0\quad  {\rm in}\quad\Omega.$$
  Our  interest is to understand what  singularity of the solution to
\begin{equation}\label{eq 1.2}
 \arraycolsep=1pt
\begin{array}{lll}
 (-\Delta)^\alpha  u+g(u)=\frac{\partial^\alpha (\omega+\delta_{x_0})}{\partial \vec{n}^\alpha}\qquad & {\rm in}\quad\bar\Omega,\\[2mm]
 \phantom{   (-\Delta)^\alpha  +g(u)}
u=0\quad & {\rm in}\quad \bar\Omega^c,
\end{array}
\end{equation}
where $x_0\in\partial\Omega$ and $\delta_{x_0}$ is the Dirac mass concentrated $x_0$ on the boundary.
Inspired by Definition \ref{weak solution GV}, it is natural to give the definition of weak solution of (\ref{eq 1.2}) as following.
\begin{definition}\label{weak solution GV1}
We say that $u$ is a weak solution of (\ref{eq 1.2}), if $u\in
L^1(\Omega)$, $g(u)\in L^1(\Omega,\rho^\alpha dx)$  and
$$
\int_\Omega [u(-\Delta)^\alpha\xi+ g(u)\xi]dx=\int_{\partial\Omega}\frac{\partial^\alpha\xi(x)}{\partial\vec{n}_x^\alpha}d\omega(x)+
\frac{\partial^\alpha\xi(x_0)}{\partial\vec{n}_{x_0}},\qquad \forall\xi\in \mathbb{X}_\alpha.
$$
\end{definition}

\begin{theorem}\label{teo 3}
Assume that   $x_0\in\partial\Omega$,  $g$ is a continuous nondecreasing function satisfying
$g(0)\ge0$,
\begin{equation}\label{g2}
 \int_1^\infty g(s)s^{-1-\frac{N+\alpha}{N-\alpha}}ds<+\infty
\end{equation}
and for some $\lambda>0$,
\begin{equation}\label{g3}
 g(s+t)\le  \lambda[g(s)+g(t)],\quad \forall s,t> 0.
\end{equation}
Then problem (\ref{eq 1.1}) admits a unique positive weak solution $v$ such that
\begin{equation}\label{eq 1.3}
  \frac{1}{c_3}\left[\rho(x)^{\alpha-1}+\frac{\rho(x)^\alpha}{|x-x_0|^{N} }\right]\le v(x)\le
  c_3\left[\rho(x)^{\alpha-1}+\frac{\rho(x)^\alpha}{|x-x_0|^{N} }\right],\qquad \forall x\in\Omega.
\end{equation}
Moreover, if assume additionally that $g$ is $C^{\beta}$ locally in $\R$ with $\beta>0$, then
$v$ is a classical solution of (\ref{1.1}).

\end{theorem}

From Theorem \ref{teo 3},  we find out a classical solution of (\ref{1.1}) with
 explosive rate  $\rho(x)^{\alpha-1}+\frac{\rho(x)^\alpha}{|x-x_0|^{N} }$,
this answers the question 3 in the first part of the introduction.

The boundary blow-up solutions of (\ref{1.1}) could be searched for by  making use of measure type data on boundary
 and the main difficulty is to do the estimate of $\mathbb{G}_\alpha[\frac{\partial^\alpha \omega}{\partial \vec{n}^\alpha}]$
 and $g(\mathbb{G}_\alpha[\frac{\partial^\alpha \omega}{\partial \vec{n}^\alpha}])$. Especially, it is dedicate to make the estimate of
 $g(\mathbb{G}_\alpha[\frac{\partial^\alpha \omega}{\partial \vec{n}^\alpha}])$ near the boundary when
 the nonlinearity $g$ is just integral-subcritical, i.e. (\ref{g1}).

 This article is organized as
follows. In Section \S2 we present some preliminaries to the Marcinkiewicz type estimate for
$\mathbb{G}_\alpha[\frac{\partial^\alpha \omega}{\partial \vec{n}^\alpha}]$ and
present the existence and uniqueness of weak solution of (\ref{eq 1.1}) when $g$ is bounded.
 Section \S3, \S4 are devoted to prove Theorem \ref{teo 1} and Theorem \ref{teo 2}.  Finally,
 we obtain one typical  solution that blows up along the boundary with different power rate.

\setcounter{equation}{0}

\section{Preliminary }

\subsection{The Marcinkiewicz type estimate}
In order to obtain the weak solution of (\ref{eq 1.1})  with integral subcritical nonlinearity,
we have to introduce the Marcinkiewicz
space and recall some related  estimate.

\begin{definition}
Let $\Theta\subset \R^N$ be a domain and $\varpi$ be a positive
Borel measure in $\Theta$. For $\kappa>1$,
$\kappa'=\kappa/(\kappa-1)$ and $u\in L^1_{loc}(\Theta,d\mu)$, we
set
\begin{equation}\label{mod M}
\|u\|_{M^\kappa(\Theta,d\varpi)}=\inf\left\{c\in[0,\infty]:\int_E|u|d\varpi\le
c\left(\int_Ed\varpi\right)^{\frac1{\kappa'}},\ \forall E\subset \Theta,\,E\
{\rm Borel}\right\}
\end{equation}
and
\begin{equation}\label{spa M}
M^\kappa(\Theta,d\varpi)=\{u\in
L_{loc}^1(\Theta,d\varpi):\|u\|_{M^\kappa(\Theta,d\varpi)}<+\infty\}.
\end{equation}
\end{definition}

The space $M^\kappa(\Theta,d\varpi)$ is called the Marcinkiewicz space of
exponent $\kappa$, or weak $L^\kappa$-space and
$\|.\|_{M^\kappa(\Theta,d\varpi)}$ is a quasi-norm. %We observe that

\begin{proposition}\label{pr 1} \cite{BBC,CC}
Assume that $1\le q< \kappa<\infty$ and $u\in L^1_{loc}(\Theta,d\varpi)$.
Then there exists  $c_4>0$ dependent of $q,\kappa$ such that
$$\int_E |u|^q d\varpi\le c_4\|u\|_{M^\kappa(\Theta,d\varpi)}\left(\int_E d\varpi\right)^{1-q/\kappa}$$
for any Borel set $E$ of $\Theta$.
\end{proposition}

Denote by $G_\alpha$ the Green kernel of $(-\Delta)^\alpha$ in $\Omega\times\Omega$ and by $\mathbb{G}_\alpha[\cdot]$ the
Green operator defined as
$$
\mathbb{G}_\alpha[\frac{\partial^\alpha \omega}{\partial \vec{n}^\alpha}](x)=\lim_{t\to0^+}\int_{\partial\Omega} G_\alpha(x,y+t\vec{n}_y)t^{-\alpha}d\omega(y).
$$
Our purpose in this subsection is to do Marcinkiewicz type estimate for $\mathbb{G}_\alpha[\frac{\partial^\alpha \omega}{\partial \vec{n}^\alpha}]$.

\begin{lemma}\label{lm 2.2}
There exists $c_5\ge1$ such that for any $x\in \Omega$,
\begin{equation}\label{2.1}
\frac1{c_5} \rho(x) ^{\alpha-1}\le \mathbb{G}_\alpha[\frac{\partial^\alpha \omega}{\partial \vec{n}^\alpha}](x)
\le c_5 \rho(x) ^{\alpha-1}.
\end{equation}

\end{lemma}
{\bf Proof.}
Since $\partial\Omega$ is $C^2$, then there exists $t_0\in(0,\frac12)$ such that
for any $x\in\Omega_t:=\{z\in\Omega,\rho(x)<t\}$ with $t< t_0$, there exists a unique $x_\partial\in \partial\Omega$ such that
$$|x-x_\partial|=\rho(x) $$
and for $t\in(0,t_0)$ letting
$$\mathcal{C}_t=\{x\in\Omega:\ \rho(x)=t \},$$
 $\mathcal{C}_t$ is $C^2$ for $t\in(0,t_0)$ and  any Borel set $E_t$ in $\mathcal{C}_t$,  there exists
unique set $E\subset \partial \Omega$ such that for any $x_t\in E_t$, there exists a unique $x\in E$ such that
\begin{equation}\label{q 1}
  |x_t-x|=t
\end{equation}
and for any $x\in E$, there exists a unique $x_t\in E_t$ satisfying (\ref{q 1}).
Moreover, for $x\in \mathcal{C}_t$ with $t\in(0,t_0)$, there exists a unique $x_\partial\in\partial\Omega$ such that
$$x=t\vec{n}_{x_\partial}+x_\partial\quad{\rm and}\quad |x-x_\partial|=t=\rho(x).$$
 Denotes by $\omega_t$ a  measure on $\mathcal{C}_t$ generated by $\omega$ such that for $t\in(0,t_0)$,
\begin{equation}\label{q 2.5}
 \omega_t(E_t)=\omega(E)\qquad {\rm for\ any\ Borel\ set\ }\ E_t\subset\mathcal{C}_t.
\end{equation}

 By compactness we only have to prove that  (\ref{2.1}) holds in a neighborhood of any point $\bar x\in\partial\Omega$ and  without loss of generality, we may assume that
$$x_\partial=0\quad {\rm and}\quad \vec{n}_{x_\partial}=e_N.$$

From \cite[Lemma 2.1]{CH}, there exists $c_6>0$ such that
\begin{equation}\label{2.2}
\mathbb{G}_\alpha[\frac{\partial^\alpha \omega}{\partial \vec{n}^\alpha}](x)\le \int_{\partial\Omega}\frac{c_6}{|x-y|^{N-\alpha}}d\omega(y),\qquad \forall x\in\Omega.
\end{equation}

Let $\phi:B'_{t_0}(0)\to \R$ such that
$(y',\phi(y')\in \partial\Omega$, where $B'_{t_0}(0)$ is the ball centered at origin with radium $t_0$ in $\R^{N-1}$.
We choose some $s_0\in(0,t_0)$ small enough, there exists $c_7\ge1$ such that
for any Borel set $E\subset B_{s_0}(0)\cap \partial\Omega$,
$$\frac1{c_7}|E'|\le \omega(E)\le c_7|E'|,$$
where
$$E'=\{y'\in\R^{N-1}:\quad (y',\phi(y'))\in E\}.$$
For $s_0>0$ small, there exists $c_8>0$ such that for $y=(y',y_N)\in B_{s_0}(0)\cap \partial\Omega$
$$|te_N-y|\ge c_8|te_N-(y',0)|=c_8\sqrt{t^2+|y'|^2} $$
Therefore,
\begin{eqnarray*}
 \int_{B_{s_0}(0)\cap \partial\Omega}\frac{1}{|te_N-y|^{N-\alpha}}d\omega(y) &\le & c_9\int_{B_{s_0}(0)\cap \partial\Omega}\frac{1}{(t^2+|y'|^2)^{\frac{N-\alpha}2}}d\omega(y)
 \\&\le &c_9\int_{B'_{s_0}(0)}\frac{1}{(t^2+|y'|^2)^{\frac{N-\alpha}2}}dy'
 \\&=&c_{10}\int_0^{s_0}\frac{s^{N-2}}{(t^2+s^2)^{\frac{N-\alpha}2}}ds
 \\&=&c_{10}t^{\alpha-1}\int_0^{\frac{s_0}t}\frac{s^{N-2}}{(1+s^2)^{\frac{N-\alpha}2}}ds
 \\&\le&c_{11} t^{\alpha-1},
\end{eqnarray*}
where $c_9,c_{10}>0$ and $c_{11}=\int_0^{+\infty}\frac{c_{10}}{(1+s^2)^{\frac{2-\alpha}2}}ds<+\infty$ since $\frac{2-\alpha}2>\frac12$.
For $y\in \partial \Omega \setminus B_{s_0}(0)$, there exists $c_{12}>0$ such that
$|te_N-y|\ge c_{12}s_0$, then
\begin{eqnarray*}
 \int_{\partial \Omega \setminus B_{s_0}(0)}\frac{1}{|te_N-y|^{N-\alpha}}d\omega(y) \le  c_{12}s_0^{\alpha-N} \int_{\partial \Omega \setminus  B_{s_0}(0)} d\omega(y)\le   c_{12}s_0^{\alpha-N}  \omega(\partial\Omega).
\end{eqnarray*}
Therefore, for $t\in(0,t_0)$,
$$\int_{\partial\Omega}\frac{1}{|te_N-y|^{N-\alpha}}d\omega(y)\le c_{12} t^{\alpha-1}.$$

For $x\in \Omega\setminus\Omega_{t_0}$ and $y\in\partial\Omega$, we observe that $ |x-y|\ge t_0$, then
$\int_{\partial\Omega}\frac{1}{|x-y|^{N-\alpha}}d\omega(y)$ is bounded by some constant dependent of $t_0$ and the diameter of $\Omega$,
thus, (\ref{2.1}) holds.

 We now prove that for $t\in(0,t_0)$,
 \begin{equation}\label{4.3-2}
\mathbb{G}_\alpha[\frac{\partial^\alpha \omega}{\partial \vec{n}^\alpha}](te_N)\ge \frac1{c_6} t^{\alpha-1}.
\end{equation}
For all $s\in(0,\frac t8)$, we have that
$$|te_N-y|> \frac t2\quad {\rm for}\quad y\in \mathcal{C}_s\cap B_{\frac t4}(se_N).$$
Therefore,
$$|te_N-y|>\frac t2 =\frac12\max\{\rho(y),\rho(te_N)\}$$
and apply \cite[Theorem 1.2]{CS} to derive that there exists $c_{13}>0$ such that for all $s\in(0,\frac t8)$
\begin{equation}\label{11.04.2}
G_\alpha(te_N,y) \ge  c_{13}\frac{\rho^\alpha (y)\rho^\alpha (te_N)}{|te_N-y|^{N}}=c_{13}\frac{t^\alpha s^\alpha}{|te_N-y|^{N}},\quad y\in \mathcal{C}_s\cap B_{\frac t4}(se_N).
\end{equation}
Thus,
$$\mathbb{G}_\alpha[s^{-\alpha}\omega_s](te_N)\ge  c_{13}\int_{\mathcal{C}_s\cap B_{\frac t4}(se_N)}\frac{ t^{\alpha}}{|te_N-y|^{N}}d\omega_s(y).$$
Denote
$$ D_{t,s}=\mathcal{C}_s\cap B_{\frac t4}(se_N)\quad{\rm and}\quad D_t=\partial\Omega\cap  B_{\frac t4}(0).$$
We observe that
$$  |te_N-y|\le c_{14}t ,\quad \forall y\in D_{t,s}$$
and
$$\frac1{c_{14}}t^{N-1}\le \omega_s(D_{t,s})\le c_{14}t^{N-1},$$
where  $c_{14}>1$,
then for any $s\in(0,\frac t8)$
\begin{eqnarray*}
\int_{\mathcal{C}_s}\frac{ t^{\alpha}}{|te_N-y|^{N}}d\omega_s(y)&\ge& \int_{D_{t,s}}\frac{ t^{\alpha}}{|x-y|^{N}}d\omega_{\mathcal{C}_s}(y)   \\[2mm]
   &\ge &  c_{15} t ^{\alpha}t^{-N}  \omega_s(D_{t,s})
     \ge c_{16} t ^{\alpha-1},
\end{eqnarray*}
which implies (\ref{4.3-2}) by passing the limit of $s\to0^+$.\qquad $\Box$

\begin{proposition}\label{pr 2.1}
Let $\mathbb{G}_\alpha[\frac{\partial^\alpha \omega}{\partial \vec{n}^\alpha}]$ given by (\ref{2.2}) and $p^*=\frac{1+\alpha}{1-\alpha}$.
 Then there exists $c_{17}>0$  such that
\begin{equation}\label{annex 0.1}
\|\mathbb{G}_\alpha[\frac{\partial^\alpha \omega}{\partial \vec{n}^\alpha}]\|_{M^{p^*}(\Omega,\rho^\alpha dx)}\le c_{17}.
\end{equation}

\end{proposition}
%%%%%%%%%%%%PROOF%%%%%%%%%%%%%%%%%%%%%%%%%%%%%%%%%%%%%%%%%%%%%%%%%%%%%%%%%%%%%%%%%%%%%%%%%%%%%%%%%%%%%%%%%%%%%%%%%%%%%%%%%%%%%%%%%%%%%%%%%%%%%
{\bf Proof.}
For any Borel set $E$ of $\Omega$ satisfying
$$0<|E|< |\Omega_{t_0}|,$$
where $\Omega_{r}=\{x\in\Omega,\ \rho(x) < r\}$ for $r>0$,
there exists $t\in(0,t_0)$ such that
$$|E|=|\Omega_{t}|. $$
Then there exists $c_{18}>0$ such that
$$ |\Omega_{t}|=\int_0^t\omega_t(\mathcal{C}_t)dt \le  c_{18}t.$$

We observe that
$$|E\setminus\Omega_{t}|=|E|-|E\cap \Omega_{t}|=|\Omega_{t}\setminus E|$$
and
$$\rho(y)\le  \rho(z),\quad \forall y\in E\setminus\Omega_{t},\ \forall z\in \Omega_{t}\setminus E. $$
Then
\begin{eqnarray*}
 \int_{E}\rho^\alpha dx&\ge &  \int_{\Omega_t} \rho^\alpha(x)dx
   =  c_{19}\int_0^t\int_{\mathcal{C}_s} s^{ \alpha}d\omega_s ds
    \ge  c_{20} t^{\alpha+1}
\end{eqnarray*}
and together with (\ref{2.1}), we deduce that
\begin{eqnarray*}
 \int_{E}\mathbb{G}_\alpha[\frac{\partial^\alpha \omega}{\partial \vec{n}^\alpha}](x)\rho^\alpha(x)dx%&\le &\int_{E} \int_{\partial\Omega}\frac{c_6}{|x-y|^{N-\alpha}}d\omega(y)\rho^\alpha(x)dx \\
   &\le & c_{21} \int_{E} \rho ^{2\alpha-1}(x) dx  \le  c_{21}\int_{\Omega_t}\rho ^{2\alpha-1}(x) dx \\
   &= & c_{22}\int_0^t\int_{\mathcal{C}_s} s^{2\alpha-1}d\omega_s ds =  c_{23} t^{2\alpha}
   \\&\le & c_{23} (\int_{E}\rho^\alpha dx)^{\frac{2\alpha}{1+\alpha}},
\end{eqnarray*}
where $c_{22},c_{23}>0$.
Therefore,
$$
\int_{E}\mathbb{G}_\alpha[\frac{\partial^\alpha \omega}{\partial \vec{n}^\alpha}](x)\rho^\alpha(x)dx\le c_{23}(\int_{E}\rho^\alpha dx)^{\frac{2\alpha}{1+\alpha}}. $$
Together with
$$ \frac{2\alpha}{1+\alpha} =\frac{p^*-1}{p^*},$$
we derive (\ref{annex 0.1}).
This completes
the proof. \qquad$\Box$

\subsection{Existence for bounded nonlinearity}

We extend Hausdorff measure $\omega$ to $\bar\Omega$ by zero inside $\Omega$, still denoting $\omega$.
For bounded $C^2$ domain, it follows \cite[p 57]{MP} that $\omega$ is a Radon measure in $\bar\Omega$.
In the approximating to weak solution of (\ref{eq 1.1}),  we consider a sequence $\{g_n\}$ of $C^1$ nonnegative  functions defined on $\R_+$
such that $g_n(0)=g(0)$,
\begin{equation}\label{2.40}
  g_n\le g_{n+1}\le g,\quad \sup_{s\in\R_+}g_n(s)=n\quad{\rm and}\quad \lim_{n\to\infty}\norm{g_n-g}_{L^\infty_{loc}(\R_+)}=0.
\end{equation}

\begin{proposition}\label{pr 2.2}
Assume that $\{g_n\}_n$ is given by (\ref{2.40}). Then
\begin{equation}\label{eq 2.1}
 \arraycolsep=1pt
\begin{array}{lll}
 (-\Delta)^\alpha  u+g_n(u)=k\frac{\partial^\alpha \omega}{\partial \vec{n}^\alpha}\quad & {\rm in}\quad\Omega,\\[2mm]
 \phantom{   (-\Delta)^\alpha  +g_n(u)}
u=0\quad & {\rm in}\quad \Omega^c
\end{array}
\end{equation}
admits a unique positive weak solution $u_{k,n}$ satisfying

$(i)$ the mapping $k\to u_{k,n}$ is increasing, the mapping $n\to u_{k,n}$ is decreasing
\begin{equation}\label{2.5}
k\mathbb{G}_\alpha[\frac{\partial^\alpha\omega}{\partial\vec{n}^\alpha}](x)-k\mathbb{G}_\alpha[g_n(k\mathbb{G}_\alpha
[\frac{\partial^\alpha\omega}{\partial\vec{n}^\alpha}])](x)\le u_{k,n}(x)\le k\mathbb{G}_\alpha[\frac{\partial^\alpha\omega}{\partial\vec{n}^\alpha}](x),\quad\forall x\in\Omega;
\end{equation}

$(ii)$
$u_{k,n}$ is a classical solution of
\begin{equation}\label{eq 2.2}
\begin{array}{lll}
\ (-\Delta)^\alpha   u+g_n(u)=0\quad & {\rm in}\quad  \Omega,\\[2mm]
\phantom{(-\Delta)^\alpha   +g_n(u)}
\ u=0\quad & {\rm in}\quad  \R^N\setminus\bar\Omega,\\[2mm]
 \phantom{}
\lim_{x\in\Omega,x\to\partial\Omega}u(x)=+\infty.
\end{array}
\end{equation}

\end{proposition}
{\bf Proof.}  Since $\omega$ is a Radon measure in $\bar\Omega$, we could apply \cite[Theorem 1.1]{CH} to
obtain that problem (\ref{eq 2.1}) admits a unique weak solution $u_{k,n}$ satisfying that $(i)$ and
$u_{k,n}$ is a classical solution of
$$
 (-\Delta)^\alpha   u+g_n(u)=0\quad  {\rm in}\quad  \Omega,\qquad \ u=0\quad  {\rm in}\quad  \R^N\setminus\bar\Omega.
 $$
 From Lemma \ref{lm 2.2} and (\ref{2.5}), there exists $c_{24}\ge1$ such that
\begin{equation}\label{eq 2.4}
\frac1{c_{24}}\rho(x)^{\alpha-1}\le u_{k,n}(x)\le c_{24}\rho(x)^{\alpha-1},\quad x\in\Omega.
\end{equation}
Therefore, $u_{k,n}$ is a classical solution of (\ref{eq 2.2}).\qquad$\Box$

\smallskip

In particular, let $g_0\equiv 0$, we have that
\begin{corollary}\label{cr 1}
$\mathbb{G}_\alpha[\frac{\partial^\alpha\omega}{\partial\vec{n}^\alpha}]$ is a classical solution of
\begin{equation}\label{eq 2.3}
\begin{array}{lll}
\qquad\quad\ \ (-\Delta)^\alpha   u=0\quad & {\rm in}\quad  \Omega,\\[2mm]
\phantom{(-\Delta)^\alpha  }
\qquad\quad\ \ u=0\quad & {\rm in}\quad  \R^N\setminus\bar\Omega,\\[2mm]
 \phantom{}
\lim_{x\in\Omega,x\to\partial\Omega}u(x)=+\infty.
\end{array}
\end{equation}
\end{corollary}

With the help of Corollary \ref{cr 1}, we are in the position to prove Proposition \ref{pr 00}.\smallskip

\noindent{\bf Proof of Proposition \ref{pr 00}.}
We first prove that $\tau_0(\alpha)\le\alpha-1$. Inversely, if $\tau_0(\alpha)>\alpha-1$, then we have that
$$1-\frac{2\alpha}{\tau_0(\alpha)}>\frac{1+\alpha}{1-\alpha}>1+2\alpha.$$
On the one hand,  it follows by \cite[Theorem 1.1]{CFQ} that for $p=\frac{1+\alpha}{1-\alpha}$, problem
\begin{equation}\label{eq 2.5}
\begin{array}{lll}
 (-\Delta)^\alpha   u+|u|^{p-1}u=0\quad & {\rm in}\quad  \Omega,\\[2mm]
\phantom{(-\Delta)^\alpha   +|u|^{p-1}u}
 u=0\quad & {\rm in}\quad  \R^N\setminus\bar\Omega,\\[2mm]
 \phantom{\  }
\lim_{x\in\Omega,x\to\partial\Omega}u(x)=+\infty
\end{array}
\end{equation}
admits a solution $w$ such that
$$\frac1{c_{25}} \le w(x)\rho^{1-\alpha}(x)\le c_{25},\quad x\in\Omega,$$
where $c_{25}>1$ and $$-\frac{2\alpha}{\frac{1+\alpha}{1-\alpha}-1}=\alpha-1.$$

On the other hand, form Corollary \ref{cr 1} we know that for any $\mu>0$, $\mu\mathbb{G}_\alpha[\frac{\partial^\alpha\omega}{\partial\vec{n}^\alpha}]$
is a super solution of problem (\ref{eq 2.5}). Furthermore, from Lemma \ref{lm 2.2},
$$\frac1{c_5}\le \mathbb{G}_\alpha[\frac{\partial^\alpha\omega}{\partial\vec{n}^\alpha}](x)\rho^{1-\alpha}\le c_5. $$
Now choosing $\mu_1=\frac1{2c_{25}c_5}$ and $\mu_2=2c_{25}c_5$, we derive that
$$\limsup_{x\in\Omega,\
x\to\partial\Omega}\mu_1\mathbb{G}_\alpha[\frac{\partial^\alpha\omega}{\partial\vec{n}^\alpha}](x)\rho(x)^{1-\alpha}<\frac1{c_{25}},\ \ \liminf_{x\in\Omega,\
x\to\partial\Omega}\mu_2\mathbb{G}_\alpha[\frac{\partial^\alpha\omega}{\partial\vec{n}^\alpha}](x)\rho(x)^{1-\alpha}> c_{25}. $$
thus, from \cite[Proposition 6.1]{CFQ},
there is no solution $u$ such that

$$\frac1{c_{25}}\le \liminf_{x\in\Omega,\
x\to\partial\Omega}u(x)\rho(x)^{1-\alpha}
\le\limsup_{x\in\Omega,\
x\to\partial\Omega}u(x)\rho(x)^{1-\alpha}\le c_{25}.
$$
The contradiction is obvious.

We finally prove that $\tau_0(\alpha)\ge\alpha-1$. Inversely, if $\tau_0(\alpha)<\alpha-1$, then we have that
$$1-\frac{2\alpha}{\tau_0(\alpha)}<\frac{1+\alpha}{1-\alpha}.$$
From \cite[Theorem 1.1]{CFQ} nonexistence $(ii)$, there is no solution $u$ of problem (\ref{eq 2.5}) with
\begin{equation}\label{p 1}
 \max\{1-\frac{2\alpha}{\tau_0(\alpha)},\frac{2\alpha}{1-\alpha}\} <p< \frac{1+\alpha}{1-\alpha}
\end{equation}
 such that
\begin{eqnarray}\label{4.0.1}
0< \liminf_{x\in\Omega,\
x\to\partial\Omega}u(x)\rho(x)^{1-\alpha}
\le\limsup_{x\in\Omega,\
x\to\partial\Omega}u(x)\rho(x)^{1-\alpha}<+\infty.
\end{eqnarray}

Let $\bar \tau=2\alpha-(1-\alpha)p$, then
$$\bar \tau>2\alpha-(1-\alpha)\frac{1+\alpha}{1-\alpha}=\alpha-1$$
and
$$\bar \tau<2\alpha-(1-\alpha)\frac{2\alpha}{1-\alpha}=0.$$
For $t_0>0$ small,   $\Omega_{t_0}=\{x\in \Omega,\rho(x)<t_0\}$ is $C^2$ and
define
\begin{equation}\label{3.3.1}
V_1(x)=\left\{ \arraycolsep=1pt
\begin{array}{lll}
 d(x)^{\bar\tau},\ & x\in \Omega_{t_0},\\[2mm]
 l(x),\ \ \ \ & x\in \Omega\setminus \Omega_{t_0},\\[2mm]
0,\ &x\in\Omega^c,
\end{array}
\right.
\end{equation}
 where   the function $l$  is  positive  such that $V_1$ is $C^2$ in $\Omega$.
 From \cite[Proposition 3.2 $(ii)$]{CFQ}, there exists $\delta_1\in(0,t_0]$ and $c_{26}>1$ such that
 $$\frac1{c_{26}}\rho(x)^{\bar\tau-2\alpha }\leq(-\Delta)^{\alpha}V_1(x)\leq
c_{26}\rho(x)^{\bar\tau-2\alpha },\ \  \forall x\in \Omega_{\delta_1}.$$
We observe that $\mathbb{G}_\alpha[\frac{\partial^\alpha\omega}{\partial\vec{n}^\alpha}]$ is a super solution of (\ref{eq 2.5})
with $p$ in (\ref{p 1}).
Now we define
$$
W_{\mu}(x)=\mathbb{G}_\alpha[\frac{\partial^\alpha\omega}{\partial\vec{n}^\alpha}]-\mu
V_{\bar \tau}(x)-\mu^2 \mathbb{G}_\alpha[1],
$$
where $\mathbb{G}_\alpha[1]$ is the solution of
$$
\begin{array}{lll}
 (-\Delta)^\alpha  u=1\quad & {\rm in}\quad  \Omega,\\[2mm]
\phantom{(-\Delta)^\alpha }
 u=0\quad & {\rm in}\quad  \R^N\setminus \Omega.
\end{array}
$$
We see that  $\bar\tau-2\alpha=(\alpha-1)p$, there exists $\mu_1>0$ such that  for $\mu\ge \mu_1$ and $x\in\Omega_{\delta_1}$,
\begin{eqnarray*}(-\Delta)^\alpha
 W_{\mu}(x)+ |W_{\mu}|^{p-1}W_{\mu}(x)&\leq&
-c_{26}\mu\rho(x)^{\bar\tau-2\alpha} +c_5^p\rho(x)^{(\alpha-1)p}\le0
\end{eqnarray*}
and
there exists $\mu_2>0$ such that  for $\mu\ge \mu_2$ and $x\in\Omega\setminus\Omega_{\delta_1}$,
\begin{eqnarray*}(-\Delta)^\alpha
 W_{\mu}(x)+ |W_{\mu}|^{p-1}W_{\mu}(x)&\leq&
c_{26}\mu \max_{\Omega\setminus\Omega_{\delta_1}}|(-\Delta)^\alpha V_1|-\mu^2 + \left(\min_{\Omega\setminus\Omega_{\delta_1}}
\mathbb{G}_\alpha[\frac{\partial^\alpha\omega}{\partial\vec{n}^\alpha}]\right)^p\le0  ,
\end{eqnarray*}
Therefore,  for $\mu=\max\{\mu_1,\mu_2\}$,  $W_{\mu}$ is a sub solution of   (\ref{eq 2.5})
with $p$ in (\ref{p 1}) and
$$c_5\rho^{\alpha-1}\ge \mathbb{G}_\alpha[\frac{\partial^\alpha\omega}{\partial\vec{n}^\alpha}]\ge W_\mu \quad {\rm and}\quad  \liminf_{x\in\Omega,\
x\to\partial\Omega} W_{\mu}(x)\rho^{1-\alpha}(x)\ge \frac1{c_5}.$$
By \cite[Theorem 2.6]{CFQ}, there exists a solution $u$   of   (\ref{eq 2.5}) with $p$ in (\ref{p 1}) satisfying
(\ref{4.0.1}). A contradiction is obtained and   the proof is complete.\hfill $\Box$

\setcounter{equation}{0}
\section{Proof of Theorem \ref{teo 1} }

\begin{lemma}\label{lm 2.3}
$(i)$ Assume that $g$ is a continuous nondecreasing function satisfying
$g(0)\ge0$ and (\ref{g1}).
Then
\begin{equation}\label{4.4}
\lim_{\rho(x)\to0^+}\mathbb{G}_\alpha[g(\mathbb{G}_\alpha[k\frac{\partial^\alpha \omega}{\partial \vec{n}^\alpha}])](x)\rho(x)^{1-\alpha}=0.
\end{equation}

$(ii)$  Assume that $p\in(0,\frac{1+\alpha}{1-\alpha})$, then there exists $c_{27}>0$ such that for any $x\in \Omega_{t}$ with $t\in(0,t_0)$,
\begin{equation}\label{4.6}
\mathbb{G}_\alpha[(\mathbb{G}_\alpha[\frac{\partial^\alpha \omega}{\partial \vec{n}^\alpha}])^p](x)
\le  c_{27} \rho(x)^{2\alpha-(1-\alpha)p} +c_{33}%\left\{\arraycolsep=1pt
%\begin{array}{lll}
%c_{32} \rho(x)^{2\alpha-(1-\alpha)p}\quad  &{\rm if}\quad p\in (\frac{2\alpha}{1-\alpha}, \frac{1+\alpha}{1-\alpha}),\\[1.5mm]
%-c_{32} \ln \rho(x) \quad  &{\rm if}\quad p= \frac{2\alpha}{1-\alpha},\\[1.5mm]
%c_{32}\quad  &{\rm if}\quad p\in(0,\frac{2\alpha}{1-\alpha}).
%\end{array}
% \right.
\end{equation}
\end{lemma}
{\bf Proof.} $(i)$
Without loss of generality, we may assume that
$$0\in\partial\Omega,\quad \vec{n}_{0}=e_N,\quad x_s=se_N$$
and we just need  prove (\ref{4.6}) and (\ref{4.4}) for $x_s$ with $s\in(0,t_0)$. It follows by Lemma \ref{lm 2.2} that
\begin{equation}\label{4.3-1}
\frac1{c_5}\rho(x)^{\alpha-1} \le \mathbb{G}_\alpha[\frac{\partial^\alpha \omega}{\partial \vec{n}^\alpha}](x)\le  c_5\rho(x)^{\alpha-1}, \qquad \forall x\in\Omega.
\end{equation}
 Combining with monotonicity of $g$, we have that
\begin{eqnarray*}
  \mathbb{G}_\alpha[g(\mathbb{G}_\alpha[k\frac{\partial^\alpha \omega}{\partial \vec{n}^\alpha}])](x_s)s^{1-\alpha} &\le& \int_\Omega G_\alpha(x_s,z)g( c_5k \rho(z)^{\alpha-1})dz s^{1-\alpha} \\
   &\le &  \int_\Omega\frac{c_5 \rho^\alpha(z) }{|x_s-z|^{N-\alpha}} g( c_5k \rho(z)^{\alpha-1})dz s^{1-\alpha} \\
   &=&  c_5  \left[\int_{ B_{ t_0}(0)\cap \Omega}\frac{s^{1-\alpha}\rho^\alpha(z)}{|x_s-z|^{N-\alpha}} g( c_5k \rho(z)^{\alpha-1})dz\right.
\\&& \left.+ \int_{ \Omega\setminus B_{ t_0 }(0)}\frac{s^{1-\alpha}\rho^\alpha(z)}{|x_s-z|^{N-\alpha}} g( c_5k \rho(z)^{\alpha-1})dz\right]
\\&:=&A_1(s)+A_2(s).
\end{eqnarray*}

Let $B'_{t_0}(0)$ be the ball with radium $t_0$ and centered at the origin in $\R^{N-1}$
and since $\partial\Omega$ is $C^2$,
there exists a $C^2$ function $\psi:B'_{\eta}(0)\to \R$ such that
$$(z',\psi(z')\in \partial\Omega\quad{\rm for\ any}\  z'\in B'_{\eta}(0),$$
where $\eta>0$.
Denote
$$\Psi(y)=(y',\psi(y'))+y_N\vec{n}_{(y',\psi(y'))},\quad \forall y=(y',y_N)\in Q_{\eta},$$
where
$$Q_{\eta}=\{z=(z',z_N)\in \R^{N-1}\times\R,\ |z'|<\eta,\ 0<z_N<\eta\}.$$
Thus, $\Psi$ is a $C^2$ diffeomorphism mapping such that
$$\Psi(y)=y,\quad \forall y=te_N,\ t\in(0,\eta).$$
Therefore,  if $t_0>0$ is chosen small enough, we have that  $\Omega\cap B_{t_0}(0)\subset \Psi(Q_{\eta})$ and
there exists $c_{28}>1$ such that for $z=\Psi(y)\in \Psi(Q_{\eta})$,
\begin{equation}\label{3.1}
\rho(z)=y_N\quad{\rm and}\quad \frac1{c_{28}} |se_N-y|\le |se_N-z|\le c_{28}|se_N-y|.
\end{equation}
Then we have that
\begin{eqnarray*}
 A_1(s) &=& \int_{ B_{ t_0}(0)\cap \Omega}\frac{s^{1-\alpha}\rho^\alpha(z)}{|x_s-z|^{N-\alpha}} g( c_5k \rho(z)^{\alpha-1})dz \\
   &\le & c_{28}\int_{ B_{ t_0}(0)\cap \Omega}\frac{s^{1-\alpha}\rho^\alpha(z)}{|x_s-\Psi^{-1}(z)|^{N-\alpha}} g( c_5k \rho(z)^{\alpha-1})dz  \\
   &\le & c_{28}\int_{Q_\eta}\frac{s^{1-\alpha}y_N^\alpha}{|x_s-y|^{N-\alpha}} g( c_5k y_N^{\alpha-1})dy.
\end{eqnarray*}

For $s\in(0, \frac18 \eta)$, we decompose $Q_{\eta}$ as following
$$
Q_{i,0}=\{z=(z',z_N)\in Q_{t_0}:\ \frac{is}2\le |z'|< \frac{(i+1)s}2,\  0< z_N<\frac s2\}
$$
and
$$
Q_{i,j}=\{z=(z',z_N)\in Q_{t_0}:\ \frac{is}2\le |z'|< \frac{(i+1)s}2,\ (j+\frac12)s\le z_N < (j+\frac32)s\},$$
where $i=0,1,\cdots N_s$, $j=1,\cdots N_s$ and $N_s$ is the largest integer number such that $N_s\le \frac{ \eta }{8s}$.

For $y\in Q_{0,1}$, we have that $\frac s2\le y_N\le \frac{3s}2$ and
\begin{eqnarray}
\int_{Q_{0,1}}\frac{s^{1-\alpha}y_N^\alpha}{|x_s-y|^{N-\alpha}} g( c_5k y_N^{\alpha-1})dy&\le&c_5s^{1+\alpha}g\ ( c_{29}k s^{\alpha-1})\int_{
B_{2/3}(e_N)} \frac{1}{|e_N-z|^{N-\alpha}}dz\nonumber
\\&=&c_5r^{-\frac{1+\alpha}{1-\alpha}}g (c_{29}kr )\int_{
B_{2/3}(e_N)} \frac{1}{|\vec{n}_0-z|^{N-\alpha}}dz\label{3.2}
\\&\le& c_{30}r^{-\frac{1+\alpha}{1-\alpha}}g (c_{29}kr ),\nonumber
\end{eqnarray}
where $r=s^{\alpha-1}$ and $c_{29} , c_{30}>0$.

For $y\in  Q_{i,0}$ with $i=0,\cdots,N_s$, we obtain that $|x_s-y|\ge \frac{i+1}{4}s$  and
\begin{eqnarray}
\int_{Q_{i,0}}\frac{s^{1-\alpha}y_N^\alpha}{|x_s-y|^{N-\alpha}} g( c_5k y_N^{\alpha-1})dy&\le& \frac{c_{31} s^{1-N}}{(1+i)^{N-\alpha}}\int_{ Q_{i,0}}
y_N^\alpha g( c_5k y_N^{\alpha-1})dy \nonumber
\\&=&\frac{c_{32}}{(1+i)^{N-\alpha}}\int_0^{\frac s2} t^\alpha g( c_5kt^{\alpha-1}) dt\label{3.3}
\\&=& \frac{c_{32}(1-\alpha)^{-1}}{(1+i)^{N-\alpha}}\int_{s^{-\frac1{N-\alpha}}}^\infty \tau^{-1-\frac{1+\alpha}{1-\alpha}} g\left( c_5k\tau\right) d\tau,\nonumber
\end{eqnarray}
where $r=s^{\alpha-1}$ and $c_{31}, c_{32}>0$.

For $y\in  Q_{i,j}$ with $i=0,\cdots,N_s$, $j=1,\cdots,N_s$  and $(i,j)\not=(0,1)$,  we derive that
$|x_s-y|\ge \frac{i+j}{4}s$  and $(j+\frac12)s\le y_N < (j+\frac32)s$
\begin{eqnarray}
\int_{Q_{i,j}}\frac{s^{1-\alpha}y_N^\alpha}{|x_s-y|^{N-\alpha}} g( c_5k y_N^{\alpha-1})dy&\le& \frac{c_{33}s^{1-N}}{(i+j)^{N-\alpha}} (js)^\alpha g\left( c_5k(js)^{\alpha-1}\right)|Q_{i,j}|\nonumber
\\&=& \frac{c_{34}j^\alpha}{(i+j)^{N-\alpha}} r^{-\frac{1+\alpha}{1-\alpha}}g (c_5kj^{\alpha-1}r ),\label{3.4}
\end{eqnarray}
where $r=s^{\alpha-1}$ and $c_{33}, c_{34}>0$.

Therefore, there exist $c_{35},c_{36}>0$ such that
$$
 A_1(s)\le  \sum_{i=0,j=0}^{N_s} \frac{c_{35}j^\alpha r^{-\frac{1+\alpha}{1-\alpha}}}{(i+j+1)^{N-\alpha}}g (c_{29}kj^{\alpha-1}r) + \sum_{i=0 }^{N_s}\frac{c_{36}}{(1+i)^{N-\alpha}}\int_{s^{-\frac1{N-\alpha}}}^\infty \tau^{-1-\frac{1+\alpha}{1-\alpha}} g\left( c_5k\tau\right) d\tau .
$$
Since
$$  \sum_{i=0}^{N_s}\frac{c_{36}}{(1+i)^{N-\alpha}}\int_{s^{-\frac1{N-\alpha}}}^\infty \tau^{-1-\frac{1+\alpha}{1-\alpha}} g\left( c_5k\tau\right) d\tau
   \le c_{37} \int_{s^{-\frac1{N-\alpha}}}^\infty \tau^{-1-\frac{1+\alpha}{1-\alpha}} g\left( c_5k\tau\right) d\tau $$
which tends to 0 as $s\to0^+$ by hypothesis (\ref{g1}).

For any $\epsilon>0$, there exists $n_\epsilon>1$ such that
$$ \sum_{i,j= n_\epsilon}^{N_s} \frac{c_{35} }{(i+j+1)^{N-\alpha}}\le \epsilon,$$
and since  $\{(j^{\alpha-1}r)^{-\frac{1+\alpha}{1-\alpha}}g (c_{29}kj^{\alpha-1}r)\}$ is uniformly bounded, we imply that
\begin{eqnarray*}
 \sum_{i,j= n_\epsilon}^{N_s} \frac{c_{35}j^\alpha}{(i+j+1)^{N-\alpha}}r^{-\frac{1+\alpha}{1-\alpha}}g (c_{29}kj^{\alpha-1}r) &=&\sum_{i,j= n_\epsilon}^{N_s} \frac{c_{35} j^{-1} }{(i+j+1)^{N-\alpha}}(j^{\alpha-1}r)^{-\frac{1+\alpha}{1-\alpha}}g (c_{29}kj^{\alpha-1}r) \\
   &\le&  c_{36}\epsilon.
\end{eqnarray*}
 For $i,j< n_\epsilon$, there exists $s_\epsilon\in (0,\eta)$ such that
$$j^\alpha r^{-\frac{1+\alpha}{1-\alpha}}g (c_{29}kj^{1-\alpha}r)\le \epsilon,\quad{\rm for}\ r\ge  s_\epsilon^{\alpha-1},$$
thus, for any $\epsilon>0$, there exists $s_\epsilon$ such that for $s\in(0,s_\epsilon)$
\begin{eqnarray*}
  \sum_{i=0,j=0}^{n_\epsilon} \frac{c_{35}j^\alpha}{(i+j+1)^{N-\alpha}}r^{-\frac{1+\alpha}{1-\alpha}}g (c_{29}kj^{\alpha-1}r) &\le&\epsilon\sum_{i=0,j=0}^{N_s} \frac{c_{35} j^{-1} }{(i+j+1)^{N-\alpha}}  \\
   &\le&  c_{37}\epsilon.
\end{eqnarray*}
Then we deduce that
$$ \lim_{s\to0^+}A_1(s)=0.$$
Therefore,
$$
\lim_{s\to0^+}\mathbb{G}_\alpha[g(\mathbb{G}_\alpha[k\frac{\partial^\alpha \omega}{\partial \vec{n}^\alpha}])](x_s)s^{1-\alpha}=0.
$$
Since $|x_s-z|\ge c_{38} t_0$ for $\Omega\setminus B_{t_0}(0)$, therefore,
\begin{eqnarray}
  A_2&\le &c_{39} s^{1-\alpha}\int_{ \Omega\setminus B_{ t_0 }(0)} \rho^\alpha(z)  g( c_5k \rho(z)^{\alpha-1})dz\nonumber \\
   &\le &  c_{39} s^{1-\alpha} \omega(\partial\Omega)\int_0^{d_0}t^\alpha g( c_5k t^{\alpha-1})dt\label{3.5}
   \\&\le& c_{40} s^{1-\alpha},\nonumber
\end{eqnarray}
where $c_{39}, c_{40}>0$ and $d_0=\max_{x\in\Omega}\rho(x)$.
Thus, (\ref{4.4}) holds.

$(ii)$ When $g(s)=s^p$ with $p\in(0,\frac{1+\alpha}{1-\alpha})$, we observe that
(\ref{3.2}) becomes that
\begin{eqnarray*}
\int_{Q_{0,1}}\frac{s^{1-\alpha}y_N^\alpha}{|x_s-y|^{N-\alpha}} ( c_5 y_N^{\alpha-1})^p dy \le c_5^ps^{1+\alpha+(\alpha-1)p},
\end{eqnarray*}
(\ref{3.3}) turns out to
\begin{eqnarray*}
\int_{Q_{i,0}}\frac{s^{1-\alpha}y_N^\alpha}{|x_s-y|^{N-\alpha}} ( c_5 y_N^{\alpha-1})^pdy &=&\frac{c_{41}}{(1+i)^{N-\alpha}}\int_0^{\frac s2} t^{\alpha+(\alpha-1)p} dt
\\&\le & \frac{c_{42}}{(1+i)^{N-\alpha}} s^{1+\alpha+(\alpha-1)p}
\end{eqnarray*}
and (\ref{3.4}) becomes that
\begin{eqnarray*}
\int_{Q_{i,j}}\frac{s^{1-\alpha}y_N^\alpha}{|x_s-y|^{N-\alpha}} ( c_5 y_N^{\alpha-1})^pdy&\le&  \frac{c_{43}}{(1+i+j)^{N-\alpha}}  s^{1+\alpha+(\alpha-1)p}.
\end{eqnarray*}
Therefore, we have that
$$\int_{ B_{t_0}(0)\cap \Omega}\frac{\rho^\alpha(y)}{|x_s-y|^{N-\alpha}} \rho(y)^{(\alpha-1)p}dy\le c_{44}s^{2\alpha+(\alpha-1)p}.$$
which, combining (\ref{3.5}), implies (\ref{4.6}).\qquad$\Box$

\smallskip

\noindent{\bf Proof of Theorem \ref{teo 1}.}
{\it To prove the existence of weak solution.} Take $\{g_n\}$  a sequence of $C^1$ nondecreasing  functions defined on $\R$
satisfying  $g_n(0)=g(0)$ and
(\ref{2.40}).
By Proposition \ref{pr 2.2},  problem (\ref{eq 2.1})
admits a unique weak solution $u_{k,n}$ such that
$$0<u_{k,n}\le \mathbb{G}_{\alpha}[\frac{\partial^\alpha \omega}{\partial \vec{n}^\alpha}]\quad{\rm  in}\quad \Omega$$
and
\begin{equation}\label{2.1.1000}
\int_\Omega [u_{k,n}(-\Delta)^\alpha\xi+g_n(u_{k,n})\xi]dx=k\int_{\partial\Omega}\frac{\partial^\alpha \xi(x)}{\partial \vec{n}_x^\alpha}d\omega(x),\quad \forall\xi\in \mathbb{X}_\alpha.
\end{equation}
%By Kato's inequatlity, see Proposition 2.4 in \cite{CV1},  for $\xi\in\mathbb{X}_\alpha$, $\xi\ge0$,
%\begin{eqnarray*}
%\int_\Omega |w|(-\Delta)^\alpha \xi dx+\int_\Omega[g_n(u_1)-g_n(u_2)]{\rm sign}(w)\xi dx\le0.
%\end{eqnarray*}

For any compact set $\mathcal{K} \subset \Omega$,
 we observe  from \cite[Lemma 3.2]{CH} that   for some $\beta\in(0,\alpha)$,
$$\norm{u_{k,n}}_{C^\beta(\mathcal{K})}\le c_{45} k.$$
Therefore, up to some subsequence, there exists $u_k$ such that
$$\lim_{n\to\infty}u_{k,n}=u_k\quad{\rm \ in}\ \Omega.$$
Then $ g_n(u_{k,n})$ converge to $g(u_k)$   in $\Omega$ as $n\to\infty$.
By Proposition \ref{pr 2.1} and (3.19) in \cite{CH}, we have that
$$u_{k,n}\to u_k\ {\rm in}\ L^1(\Omega),\quad \norm{g_n(u_n)}_{L^1(\Omega,\rho^\alpha dx)}\le c_{46}\norm{\mathbb{G}_{\alpha}[\frac{\partial^\alpha |\nu|}{\partial \vec{n}^\alpha}]}_{L^1(\Omega)}$$
and
$$  m(\lambda)\leq c_{47}\lambda^{-\frac{1+\alpha}{1-\alpha}} \ \quad {\rm for}\ \ \  \lambda>\lambda_0,$$
where
$$  m(\lambda)=\int_{  S_\lambda}\rho_{\partial\Omega}^{\alpha}(x)dx
\quad{\rm  with}\quad   S_\lambda=\{x\in\Omega: \mathbb{G}_{\alpha}[\frac{\partial^\alpha |\nu|}{\partial \vec{n}^\alpha}]>\lambda\}.$$
For any Borel
set $E\subset\Omega$, we have that
$$
\displaystyle\begin{array}{lll}
\displaystyle\int_{E}|g_n(u_n)|\rho_{\partial\Omega}^\alpha(x) dx\le \int_{E\cap\tilde S^c_{\frac{\lambda}{k}}}g\left(k\mathbb{G}_{\alpha}[\frac{\partial^\alpha |\nu|}{\partial \vec{n}^\alpha}]\right)\rho_{\partial\Omega}^\alpha(x) dx+\int_{E\cap \tilde S_{\frac{\lambda}{k}}}g\left(k\mathbb{G}_{\alpha}[\frac{\partial^\alpha |\nu|}{\partial \vec{n}^\alpha}]\right)\rho_{\partial\Omega}^\alpha(x) dx
\\[4mm]\phantom{\int_{E}|g(u_t)|\rho^{\alpha}_{\partial\Omega}(x)dx}
\displaystyle\leq \tilde g\left(\frac{\lambda}{k}\right)\int_E\rho^{\alpha}_{\partial\Omega}(x)dx+\int_{\tilde S_{\frac{\lambda}{k}}}\tilde g\left(k\mathbb{G}_{\alpha}[\frac{\partial^\alpha |\nu|}{\partial \vec{n}^\alpha}]\right)\rho^{\alpha}_{\partial\Omega}(x)dx
\\[4mm]\phantom{\int_{E}|g(u_t)|\rho^{\alpha}_{\partial\Omega}(x)dx}
\displaystyle\leq \tilde
g\left(\frac{\lambda}{k}\right)\int_E\rho^{\alpha}_{\partial\Omega}(x)dx+\tilde m\left(\frac{\lambda}{k}\right) \tilde g\left(\frac{\lambda}{k}\right)+\int_{\frac{\lambda}{k}}^\infty\tilde m(s)d\tilde
g(s),
\end{array}
$$
where $\tilde g(r)=g(|r|)-g(-|r|)$.

On the other hand,
$$\int_{\frac{\lambda}{k}}^\infty \tilde g(s)d\tilde m(s)=\lim_{T\to\infty}\int_{\frac{\lambda}{k}}^T \tilde g(s)d \tilde m(s).
$$
Thus,
$$\displaystyle\begin{array}{lll}
\displaystyle \tilde m\left(\frac{\lambda}{k}\right) \tilde g\left(\frac{\lambda}{k}\right)+ \int_{\frac{\lambda}{k}}^T \tilde m(s)d\tilde g(s) \le c_{47}\tilde g\left(\frac{\lambda}{k}\right)\left(\frac{\lambda}{k}\right)^{-\frac{1+\alpha}{1-\alpha}}+c_{47}\int_{\frac{\lambda}{k}}^T s^{-\frac{1+\alpha}{1-\alpha}}d\tilde g(s)
\\[4mm]\phantom{-----\ \int_{\lambda}^T \tilde g(s)d\omega(s)}\displaystyle
\leq c_{48}T^{-\frac{1+\alpha}{1-\alpha}}\tilde
g(T)+\frac{c_{49}}{\frac{1+\alpha}{1-\alpha}+1}\int_{\frac{\lambda}{k}}^T
s^{-1-\frac{1+\alpha}{1-\alpha}}\tilde g(s)ds.
\end{array}$$
By assumption (\ref{g1}) and \cite[Lemma 3.4]{CH} with $p=\frac{1+\alpha}{1-\alpha}$,  $T^{-\frac{N+\alpha}{N-\alpha}}\tilde g(T)\to 0$ when $T\to\infty$, therefore,
$$\tilde m\left(\frac{\lambda}{k}\right) \tilde g\left(\frac{\lambda}{k}\right)+ \int_{\frac{\lambda}{k}}^\infty \tilde m(s)\ d\tilde g(s)\leq \frac{c_{49}}{\frac{1+\alpha}{1-\alpha}+1}\int_{\frac{\lambda}{k}}^\infty s^{-1-\frac{1+\alpha}{1-\alpha}}\tilde g(s)ds.
$$
Notice that the above quantity on the right-hand side tends to $0$
when $\lambda\to\infty$. The conclusion follows: for any
$\epsilon>0$ there exists $\lambda>0$ such that
$$\frac{c_{49}}{\frac{1+\alpha}{1-\alpha}+1}\int_{\frac{\lambda}{k}}^\infty s^{-1-\frac{1+\alpha}{1-\alpha}}\tilde g(s)ds\leq \frac{\epsilon}{2}.
$$
For $\lambda$ fixed,  there exists $\delta>0$ such that
$$\int_E\rho_{\partial\Omega}^\alpha(x) dx\leq \delta\Longrightarrow \tilde g\left(\frac{\lambda}{k}\right)\int_E\rho_{\partial\Omega}^\alpha(x) dx\leq\frac{\epsilon}{2},
$$
which implies that $\{g_n\circ u_n\}$ is uniformly integrable in
$L^1(\Omega,\rho_{\partial\Omega}^\alpha dx)$. Then $g_n\circ u_n\to g\circ u_\nu$ in
$L^1(\Omega,\rho_{\partial\Omega}^\alpha dx)$ by Vitali convergence theorem.

Passing to the limit as
$n\to +\infty$ in the identity (\ref{2.1.1000}),
it implies that
$$\int_\Omega [u_k(-\Delta)^\alpha\xi+g(u_k)\xi]dx=k\int_{\partial\Omega}\frac{\partial^\alpha\xi(x)}{\partial\vec{n}^\alpha_x}d\omega(x),\quad \forall\xi\in\mathbb{ X}_\alpha.$$
Then $u_k$ is a weak solution of (\ref{eq 1.1}). Moreover,
it follows by the fact
 \begin{equation}\label{2.6}
 k\mathbb{G}_\alpha[\frac{\partial^\alpha \omega}{\partial \vec{n}^\alpha}]-g(k\mathbb{G}_\alpha[\frac{\partial^\alpha \omega}{\partial \vec{n}^\alpha}])\le u_k\le  k\mathbb{G}_\alpha[\frac{\partial^\alpha \omega}{\partial \vec{n}^\alpha}]\quad {\rm in}\ \Omega.
 \end{equation}

{\it   Uniqueness of weak solution.}
 Let $u_1,u_2$  be two weak solutions of
(\ref{eq 1.1}) and $w=u_1-u_2$. Then $(-\Delta)^\alpha
w=g_n(u_2)-g_n(u_1)$ and $g_n(u_2)-g_n(u_1)\in L^1(\Omega,\rho^\alpha dx)$.
By Kato's inequatlity, see  \cite[Proposition 2.4]{CV1},  for $\xi\in\mathbb{X}_\alpha$, $\xi\ge0$, we have that
\begin{eqnarray*}
\int_\Omega |w|(-\Delta)^\alpha \xi dx+\int_\Omega[g_n(u_1)-g_n(u_2)]{\rm sign}(w)\xi dx\le0.
\end{eqnarray*}
Combining with $\int_\Omega[g_n(u_1)-g_n(u_2)]{\rm sign}(w)\xi dx\ge0$,
then we have
$$w=0\quad {\rm a.e.\ in}\ \ \Omega.$$

{\it Regularity of $u_{k,n}$ and $u_k$.} Since $g_n$ is $C^1$ in $\R$, then by \cite[Lemma 3.2]{CH}, we have
\begin{equation}\label{2.0.10}
\norm{u_{k,n}}_{C^{2\alpha+\beta}(\mathcal{K})}\le c_{50}k,
\end{equation}
for any compact set $\mathcal{K}$ and some $\beta\in(0,\alpha)$. Then  $u_{k,n}$ is $C^{2\alpha+\beta}$ locally in $\Omega$.
Together with the fact that $u_{n,k}$ is classical solution of (\ref{eq 2.2}), we derive by Theorem 2.2 in \cite{CFQ} that
$u_k$ is a classical solution of (\ref{eq 2.2}).

{\it To prove  (\ref{1.33}.) }
Plugging (\ref{4.4}) and (\ref{4.3-1}) into (\ref{2.6}), we obtain that (\ref{1.33}).   \qquad$\Box$

\setcounter{equation}{0}
\section{Proof of Theorem \ref{teo 2} }
\subsection{Strong singularity for $p\in(1+2\alpha ,\frac{1+\alpha}{1-\alpha})$}

In this subsection, we consider the limit of $\{u_k\}$ as $k\to\infty$, where $u_k$ is the weak solution of
$$
\arraycolsep=1pt
\begin{array}{lll}
 (-\Delta)^\alpha   u+u^p=k\frac{\partial^\alpha\omega}{\partial \vec{n}^\alpha}\quad  &{\rm in}\quad\ \ \bar\Omega,\\[3mm]
 \phantom{-----\ }
 u=0\quad &{\rm in}\quad\ \ \bar\Omega^c,
 \end{array}
 $$
here  $p\in(1+ 2\alpha,\frac{1+\alpha}{1-\alpha})$.
From Theorem \ref{teo 1}, we know that $k\mapsto u_k$ is increasing and $u_k$ is a classical solution of (\ref{1.3}).

In order to control the limit of $\{u_k\}$ as $k\to\infty$, we have to obtain barrier function, i.e. a suitable super solution of (\ref{1.3}). To this end, we consider $C^2$ function $w_p$ satisfying
\begin{equation}\label{4.2}
  w_p(x)=\left\{
  \begin{array}{lll}
 \rho(x)^{-\frac{2\alpha}{p-1}},\quad&{\rm for}\ \ x\in \Omega,\\[3mm]
 \phantom{}
 0,\quad & {\rm for}\ \  x\in\Omega^c.
 \end{array}
  \right.
\end{equation}
We see that $ w_p\in L^1(\Omega)$ if $\frac{2\alpha}{p-1}<1$, i.e. $p>1+2\alpha$.

\begin{lemma}\label{lm 4.1}
Assume that $p\in(1+ 2\alpha ,\frac{1+\alpha}{1-\alpha})$ and $w_p$ is defined in (\ref{4.2}).
Then there exists $\lambda_0>0$ such that $\lambda_0 w_p$ is a super solution of
(\ref{1.3}).

\end{lemma}
{\bf Proof.} For $p\in(1+ 2\alpha ,\frac{1+\alpha}{1-\alpha})$, we have that $-\frac{2\alpha}{p-1}\in (-1,0)$
and  from \cite[Proposition 3.2]{CFQ}, it shows that there exists $c(p)<0$ such that
$$(-\Delta)^\alpha w_p(x)\ge c(p)\rho(x)^{-\frac{2\alpha}{p-1}-2\alpha},\quad x\in \Omega. $$
Thus, taking $\lambda_0=|c(p)|^{\frac{1}{p-1}}$, we derive that
$$ (-\Delta)^\alpha (\lambda_0 w_p)+ (\lambda_0w_p)^p\ge0\quad{\rm in}\quad \Omega.$$
Together with $\lambda_0 w_p=0$ in $\Omega^c$, $\lambda_0 w_p$ is a super solution of
(\ref{1.3}).
The proof ends. \qquad$\Box$

\smallskip
We observe that the super solution $\lambda_0w_p$ constructed in Lemma \ref{lm 4.1} provide a upper bound for $u_\infty$.
%In order to obtain the lower bound for $u_\infty$.

\medskip

\noindent{\bf Proof of Theorem  \ref{teo 2} $(i)$.}
For $p\in(1+ 2\alpha,\frac{1+\alpha}{1-\alpha})$, we have that
$$-\frac{2\alpha}{p-1}\in(-1,-1+\alpha)$$
and it follows by (\ref{4.3-1}) that
$$u_k(x)\le k\mathbb{G}_\alpha[\frac{\partial^\alpha \omega}{\partial \vec{n}^\alpha}](x)\le \frac{c_5k}{|x|^{1-\alpha}},\quad x\in\Omega.$$
Then
$\lim_{x\in\Omega,\rho(x)\to0}\frac{u_k(x)}{ w_p(x)}=0$
and we claim that
$$u_k\le \lambda_0w_p\quad{\rm in}\quad\Omega. $$
In fact, if it fails, then there  exists $z_0\in\Omega$ such that
$$(u_k-\lambda_0w_p)(z_0)=\inf_{\Omega}(u_k-\lambda_0w_p)={\rm ess}\inf_{\R^N}(u_k-\lambda_0w_p)<0.$$
Then we have
$(-\Delta)^\alpha(u_k-\lambda_0w_p)(z_0)<0$, which contradicts
the fact that
$$(-\Delta)^\alpha(u_k-\lambda_0w_p)(z_0)=\lambda_0w_p^p(z_0)-u_k^p(z_0)>0.$$

 By monotonicity of the mapping $k\to u_k$, there holds
$$u_\infty(x):=\lim_{k\to\infty} u_k(x),\quad x\in\R^N\setminus\{0\},$$
 which is a classical solution
of (\ref{1.3}) and
$$u_\infty(x)\le \lambda_0w_p(x)=  \lambda_0 \rho(x)^{-\frac{2\alpha}{p-1}},\quad \forall x\in\Omega.$$
By applying Stability Theorem \cite[Theorem 2.4]{CFQ}, we obtain that $u_\infty$  is a classical solution of (\ref{1.3}).

Finally, we claim that there exists $c_{51}>0$  such that for $t\in(0,t_0)$,
\begin{equation}\label{13-08-0}
 u_\infty(x)\ge c_{51}t^{-\frac{2\alpha}{p-1}},\quad \forall x\in \mathcal{C}_t.
\end{equation}

Indeed, let $t_k=(\sigma^{-1} k)^{\frac{p-1}{(1-\alpha)p-1-\alpha}}$, where $\sigma>0$ will be chosen later, then $k=\sigma t_k^{\frac{(1-\alpha)p-1-\alpha}{p-1}}$
and for $x\in \mathcal{C}_t$, we apply Lemma \ref{lm 2.3} with $p\in(1+2\alpha ,\frac{1+\alpha}{1-\alpha})$ that
\begin{eqnarray*}
 u_k(x)&\ge & k\mathbb{G}_\alpha[\frac{\partial^\alpha \omega}{\partial \vec{n}^\alpha}](x)-k^p\mathbb{G}_\alpha[(\mathbb{G}_\alpha[\frac{\partial^\alpha \omega}{\partial \vec{n}^\alpha}])^p](x) \\
 &\ge& c_5kt^{\alpha-N}[1-c_{52}k^{p-1}t^{(\alpha-1)p+\alpha+1}]
  \\&\ge& c_5\sigma t_k^{-\frac{2\alpha }{p-1}}[1-c_{52}\sigma^{p-1} t_k^{p-1}(t_k/2)^{(\alpha-1)p+\alpha+1}]
\\&\ge& c_5\sigma t_k^{-\frac{2\alpha }{p-1}}[1-c_{52}\sigma^{p-1} 2^{(1-\alpha)p-\alpha-1}]
 \\&\ge& \frac{c_5\sigma}{2}\rho(x)^{-\frac{2\alpha }{p-1}},
\end{eqnarray*}
where we choose $\sigma$ such that $c_{52}\sigma^{p-1} 2^{(N-\alpha)p-\alpha-N}=\frac12$.
Then for any $x\in \Omega$, there exists $k>0$ such that
$x\in \Omega$ and then
$$u_\infty(x)\ge u_k(x)\ge \frac{c_5\sigma}{2}\rho(x)^{-\frac{2\alpha }{p-1}},\qquad \forall x\in \Omega.$$
This ends the proof. \qquad$\Box$

\subsection{The limit of $\{u_k\}$ blows up when $p\in(0,1+ 2\alpha ]$}

In this subsection, we derive the blow-up behavior of the limit of $\{u_k\}$  when $p\in(0,1+ 2\alpha ]$. To this end,
we have to do more  estimate for $u_k$.
\begin{lemma}\label{lm 3.2}
Assume  that $g(s)=s^p$ with $p\in(1,1+2\alpha]$ and $u_k$ is the solution of (\ref{eq 1.1}) obtained by Theorem \ref{teo 1}.
Then there exist $c_{52}>0$, $r_0\in(0,\frac14)$ and $\{r_k\}_k\subset(0,r_0)$ satisfying $r_k\to0$ as $k\to\infty$ such that
\begin{equation}\label{4.3.1}
u_k(x)\ge\frac{c_{52}}{\rho(x)},\qquad \forall x\in\Omega_{t_0}.
\end{equation}
\end{lemma}
{\bf Proof.} We divide $\alpha,p$ into 4 cases:\\
 Case I: $1<\frac{2\alpha}{1-\alpha}<1+2\alpha$ and $p\in [\frac{2\alpha}{1-\alpha},1+2\alpha]$;\quad
Case II: $1<\frac{2\alpha}{1-\alpha}<1+2\alpha$ and $p\in [1,\frac{2\alpha}{1-\alpha})$;\\
Case III: $\frac{2\alpha}{1-\alpha}\le  1$ and $p\in (1,1+2\alpha]$;\qquad \qquad\ \ \
Case IV: $\frac{2\alpha}{1-\alpha}\ge 1+2\alpha$ and $p\in (1,1+2\alpha]$.

{\it To prove (\ref{4.3.1}) in Case I and Case III.} Let $r_j=j^{-\frac1\alpha}$ with $j\in(k_0,k)$, then $j=r_j^{-\alpha}$.
Applying Lemma \ref{lm 2.3} $(ii)$, for $p\ge \frac{2\alpha}{1-\alpha}$, we have that for $x\in \Omega_{r_j} \setminus \Omega_{\frac{r_j}2}$,
\begin{eqnarray*}
 u_j(x)&\ge & j\mathbb{G}_\alpha[\frac{\partial^\alpha \omega}{\partial \vec{n}^\alpha}](x)-j^p\mathbb{G}_\alpha[(\mathbb{G}_\alpha[\frac{\partial^\alpha \omega}{\partial \vec{n}^\alpha}])^p](x) \\
 &\ge& c_5^{-1}jr_j^{\alpha-1}-c_{27}j^p\rho(x)^{(\alpha-1)p+2\alpha}
\\&\ge& c_5^{-1}r_j^{-1}-c_{27}r_j^{-\alpha p-(1-\alpha)p+2\alpha}
 \\&\ge& \frac{1}{2c_5}\rho(x)^{-1},
\end{eqnarray*}
where  the last inequality holds since $-\alpha p-(1-\alpha)p+2\alpha>-1$ and $r_j\to0$ as $j\to\infty$.
Then for any $x\in \Omega_{t_0} $, there exists $j\in (k_0,k)$ such that
$x\in \Omega_{t_j} \setminus \Omega_{\frac{t_j}2}$ and then
$$u_k(x)\ge u_j(x)\ge \frac{1}{2c_5}\rho^{-1}(x),\qquad \forall x\in \Omega_{t_0}.$$

{\it To prove (\ref{4.3.1}) in Case II and Case IV.} Let $r_j=j^{-\frac1\alpha}$ with $j\in(k_0,k)$, then $j=r_j^{-\alpha}$
and for $x\in \Omega_{r_j} \setminus \Omega_{\frac{r_j}2}$, we have that
\begin{eqnarray*}
 u_j(x)&\ge & j\mathbb{G}_\alpha[\frac{\partial^\alpha \omega}{\partial \vec{n}^\alpha}](x)-j^p\mathbb{G}_\alpha[(\mathbb{G}_\alpha[\frac{\partial^\alpha \omega}{\partial \vec{n}^\alpha}])^p](x) \\
 &\ge& c_5^{-1}j\rho(x)^{\alpha-1}-c_{27}j^p
\\&\ge& c_5^{-1}r_j^{-1}-c_{27}r_j^{-\alpha p}
 \\&\ge& \frac{1}{2c_5}\rho(x)^{-1},
\end{eqnarray*}
where  the last inequality holds since $-\alpha p>-1$ and $r_j\to0$ as $j\to\infty$.
For any $x\in \Omega_{t_0} $, there exists $j\in (k_0,k)$ such that
$x\in \Omega_{r_j} \setminus \Omega_{\frac{r_j}2}$ and then
$$u_k(x)\ge u_j(x)\ge \frac{1}{2c_5}\rho(x)^{-1},\qquad \forall x\in \Omega_{t_0}.$$
%
%{\it To prove (\ref{4.3.1}) in the case of $p=1+ 2\alpha $ if $\frac{2\alpha}{1-\alpha}<1+2\alpha$.} Let $s_j=j^{-\frac1\alpha}$ and $r_j=\frac{s_j}{[-\log(s_j)]^{\frac1\alpha}}$, then $j=s_j^{-\alpha}$
%and applied Lemma \ref{lm 2.3} part $(ii)$ for $x\in A_{r_0}\cap \left[B_{r_j}(0)\setminus B_{\frac{r_j}{2}}(0)\right]$,
%\begin{eqnarray*}
% u_j(x)&\ge & j\mathbb{G}_\alpha[\frac{\partial^\alpha \omega}{\partial \vec{n}^\alpha}](x)-j^p\mathbb{G}_\alpha[(\mathbb{G}_\alpha[\frac{\partial^\alpha \omega}{\partial \vec{n}^\alpha}])^p](x) \\
% &\ge& c_5^{-1}j\rho(x)^{\alpha-1}-c_{27} j^p\rho(x)^{(\alpha-1)p+2\alpha}
%\\&\ge& c_5^{-1}s_j^{-1} (-\log s_j)^{\frac{1-\alpha}{\alpha}}-c_{42} s_j^{-1} (-\log s_j)^{\frac{(1-\alpha)p-2\alpha}{\alpha}}
% \\&=& c_5^{-1}s_j^{-1} (-\log s_j)^{\frac{1-\alpha}{\alpha}}\left[1-c_{42}(-\log s_j)^{\frac{(1-\alpha)p-2\alpha}{\alpha}-\frac{1-\alpha}{\alpha}}\right]
% \\&\ge& c_5^{-1}\frac{r_j^{-1}}{-\log s_j}[1-c_{42}(-\log s_j)^{\frac{(1-\alpha)p-1-\alpha}{\alpha}}]
% \\&\ge&  c_5\frac1{\rho(x)\log (\rho(x))},
%\end{eqnarray*}
%where $c_{42}>0$ and we used the facts that $\log(s_j)\le c\log r_j\le c\log \rho(x)$ and $\frac{(1-\alpha)p-1-\alpha}{\alpha}<0$.
%Then for any $x\in \Omega_{t_0}$, there exists $j\in (k_0,k)$ such that
%$x\in \Omega_{t_j} \setminus \Omega_{\frac{t_j}2}$ and then
%$$u_k(x)\ge c_5\frac1{-\rho(x)\log (\rho(x))} ,\qquad x\in\Omega_{t_0}.$$
The proof ends. \qquad$\Box$
\smallskip

\noindent{\bf Proof of Theorem  \ref{teo 2} $(ii)$.}
It derives by Lemma \ref{lm 3.2} that
\begin{equation}\label{3.2.3}
\pi_k:=\int_{\Omega_{t_0}\setminus \Omega_{\frac{t_k}2}}u_k(x)\ge c_{52}\int_{\Omega_{t_0}\setminus \Omega_{\frac{t_k}2}}\rho^{-1}(x)dx\to\infty\quad {\rm as}\ k\to\infty.
\end{equation}
Then
\begin{equation}\label{5.1}
\arraycolsep=1pt
\begin{array}{lll}
 (-\Delta)^\alpha  u+u^p=0 \quad & {\rm in}\quad \Omega\setminus \Omega_{t_0},\\[2mm]
 \phantom{  (-\Delta)^\alpha  +u^p}
u=0  \quad & {\rm in}\quad \R^N \setminus \Omega,\\[2mm]
\phantom{ (-\Delta)^\alpha  +u^p}
u=u_k  \quad & {\rm in}\quad \Omega_{t_0}
\end{array}
\end{equation}
admits a unique solution $w_k$.
By Comparison Principle,
\begin{equation}\label{4.1.3}
 u_{k}\ge w_k\quad {\rm in}\quad B_{\varrho_0}(y_0).
\end{equation}
Let  $\tilde w_k=w_k-u_k\chi_{\Omega\setminus \Omega_{t_0}},$
then $\tilde w_k=w_k$ in $\Omega\setminus \Omega_{t_0}$ and for $x\in \Omega\setminus \Omega_{t_0}$,
$$
\arraycolsep=1pt
\begin{array}{lll}
(-\Delta)^\alpha \tilde w_k(x) =
-\lim_{\epsilon\to0^+}\int_{B_{\varrho_0}(y_0)\setminus B_\epsilon(x)}\frac{w_k(z)-w_k(x)}{|z-x|^{N+2\alpha}}dz
\\[3mm]\phantom{--------}+\lim_{\epsilon\to0^+}\int_{B_{\varrho_0}^c(y_0)\setminus B_\epsilon(x)}\frac{w_k(x)}{|z-x|^{N+2\alpha}}dz
\\[3mm]\phantom{------}
=-\lim_{\epsilon\to0^+}\int_{\R^N\setminus B_\epsilon(x)}\frac{w_k(z)-w_k(x)}{|z-x|^{N+2\alpha}}dz +\int_{B_{r_0}(0)}\frac{u_k(z)}{|z-x|^{N+2\alpha}}dz
\\[3mm]\phantom{------}\ge(-\Delta)^\alpha w_k(x)+c_{53}\pi_k,
\end{array}
$$
where $c_{53}=(|y_0|+r_0)^{-N-2\alpha}$ and the last inequality follows by the fact of $$|z-x|\le |x|+|z|\le |y_0|+r_0\quad
{\rm for}\ \forall z\in B_{\frac14}(0),\ \ \forall x\in \Omega\setminus \Omega_{t_0}.$$
Therefore,
\begin{eqnarray*}
(-\Delta)^\alpha \tilde w_k(x)+\tilde w_k^p(x) &\ge&  (-\Delta)^\alpha w_k(x)+w_k^p(x)+ c_{53}\pi_k \\
        &=&c_{53}\pi_k, \qquad\forall x\in \Omega\setminus \Omega_{t_0},
     \end{eqnarray*}
that is,  $\tilde w_k$  is a super solution of
\begin{equation}\label{4.1.2}
\arraycolsep=1pt
\begin{array}{lll}
\displaystyle (-\Delta)^\alpha  u+u^p=c_{53}\pi_k \quad & {\rm in}\quad  \Omega\setminus \Omega_{t_0},\\[2mm]
 \phantom{  (-\Delta)^\alpha  +u^{p,}}
u=0  \quad & {\rm in}\quad (\Omega\setminus \Omega_{t_0}) ^c .
\end{array}
\end{equation}
Let $\eta_1$ be the solution of
$$
\arraycolsep=1pt
\begin{array}{lll}
 (-\Delta)^\alpha  u=1 \quad & {\rm in}\quad \Omega\setminus \Omega_{t_0},\\[2mm]
 \phantom{  (-\Delta)^\alpha  }
u=0  \quad & {\rm in}\quad  (\Omega\setminus \Omega_{t_0})^c.
\end{array}
$$
Then
 $(c_{53}\pi_k)^{\frac1p} \frac{\eta_1}{2\max_{\R^N}\eta_1}$ is sub solution of (\ref{4.1.2}) for $k $ large enough. By Comparison Principle,  we have that
 $$\tilde w_k(x)\ge (c_{53}\pi_k)^{\frac1p} \frac{\eta_1(x)}{2\max_{\R^N}\eta_1},\quad \forall  x\in B_{\varrho_0}(y_0),$$
 which implies that
 $$w_k(y)\ge c_{54} (c_{53}\pi_k)^{\frac1p},\qquad \forall y\in \Omega\setminus \Omega_{t_0},$$
 where $c_{54}=\min_{x\in B_{\varrho_0}(y_0)}\frac{\eta_1(x)}{2\max_{\R^N}\eta_1}$.
 Therefore, (\ref{4.1.3}) and (\ref{3.2.3}) imply that $$\lim_{k\to\infty}u_{k}(y)\ge \lim_{k\to\infty}w_k(y)=\infty,\qquad\forall \Omega\setminus \Omega_{t_0}.$$
Similarly, we can prove
$$\lim_{k\to\infty}u_{k}(y)\ge \lim_{k\to\infty}w_k(y)=\infty,\qquad\forall y\in \Omega.$$
The proof ends.\qquad$\Box$

\setcounter{equation}{0}
\section{ Proof of Theorem \ref{teo 3}  }

In this section, we are devoted to consider the solution of (\ref{1.1})
with different blow-up speeds. Without loss of generality, we assume that
$$x_0=0\in\partial\Omega.$$

 For $k,j\in\N$, donte
\begin{equation}\label{5.2}
\nu=\omega+\delta_0,
\end{equation}
and then
$$\mathbb{G}_\alpha[\frac{\partial^\alpha \nu}{\partial \vec{n}^\alpha}]=\mathbb{G}_\alpha[\frac{\partial^\alpha \delta_0}{\partial \vec{n}^\alpha}]
+\mathbb{G}_\alpha[\frac{\partial^\alpha \omega}{\partial \vec{n}^\alpha}],$$
which  is a $\alpha-$harmonic function of (\ref{eq 2.3}), since $\mathbb{G}_\alpha[\frac{\partial^\alpha \omega}{\partial \vec{n}^\alpha}]$
is  a $\alpha-$harmonic function of (\ref{eq 2.3}) and $\mathbb{G}_\alpha[\frac{\partial^\alpha \delta_0}{\partial \vec{n}^\alpha}]$
is $\alpha-$ harmonic.

\begin{proposition}\label{pr 4.1}
Let $\mathbb{G}_\alpha[\frac{\partial^\alpha \nu}{\partial \vec{n}^\alpha}]$ given by (\ref{5.2}) and $p^*_N=\frac{N+\alpha}{N-\alpha}$.
 Then there exists $c_{55}>0$  such that
\begin{equation}\label{5.3}
\|\mathbb{G}_\alpha[\frac{\partial^\alpha \nu}{\partial \vec{n}^\alpha}]\|_{M^{p^*_N}(\Omega,\rho^\alpha dx)}\le c_{55}.
\end{equation}

\end{proposition}
{\bf Proof.} Since
\begin{equation}\label{p m0}
\|\mathbb{G}_\alpha[\frac{\partial^\alpha \nu}{\partial \vec{n}^\alpha}]\|_{M^{p^*_N}(\Omega,\rho^\alpha dx)}\le
k \|\mathbb{G}_\alpha[\frac{\partial^\alpha \delta_0}{\partial \vec{n}^\alpha}]\|_{M^{p^*_N}(\Omega,\rho^\alpha dx)}
+j\|\mathbb{G}_\alpha[\frac{\partial^\alpha \omega}{\partial \vec{n}^\alpha}]\|_{M^{p^*_N}(\Omega,\rho^\alpha dx)}.
\end{equation}
From Proposition \ref{pr 2.1} and $p*=\frac{1+\alpha}{1-\alpha}>\frac{N+\alpha}{N-\alpha}$, on the one hand, we have that
\begin{eqnarray}
 \|\mathbb{G}_\alpha[\frac{\partial^\alpha \delta_0}{\partial \vec{n}^\alpha}]\|_{M^{p^*_N}(\Omega,\rho^\alpha dx)}&\le &
  c_{56}\|\mathbb{G}_\alpha[\frac{\partial^\alpha \delta_0}{\partial \vec{n}^\alpha}]\|_{M^{p^*}(\Omega,\rho^\alpha dx)} \nonumber\\
&\le & c_{56}c_{55},\label{5.4}
\end{eqnarray}
On the other hand, for $\lambda>0$, denote
$$
   S_\lambda=\{x\in \Omega:\mathbb{G}_{\alpha}[\frac{\partial^\alpha \delta_0}{\partial \vec{n}^\alpha}](x)>\lambda\}\quad{\rm and}\quad
    m(\lambda)=\int_{S_\lambda} \rho_{\partial\Omega}^\alpha(x)dx
$$
and from \cite[Lemma 3.3]{CH} with $\nu=\delta_0$,  there exist $\lambda_0>1$ and $c_{57}>0$ such that for any $\lambda\ge \lambda_0$,
$$
\tilde m(\lambda)\le c_{57}\lambda^{-\frac{N+\alpha}{N-\alpha}}.
$$
which implies that
\begin{equation}\label{5.5}
 \|\mathbb{G}_\alpha[\frac{\partial^\alpha \delta_0}{\partial \vec{n}^\alpha}]\|_{M^{p^*_N}(\Omega,\rho^\alpha dx)}\le c_{55},
\end{equation}
Thus, (\ref{5.3}) follows by (\ref{5.4}) and (\ref{5.5}).
Combining \cite[Theorem 1.1]{CH} and the proof of Proposition \ref{pr 2.2}, we have following result.
\begin{proposition}\label{pr 4.2}
Assume that $\{g_n\}_n$ is given by (\ref{2.40}). Then
\begin{equation}\label{q 2.1}
 \arraycolsep=1pt
\begin{array}{lll}
 (-\Delta)^\alpha  u+g_n(u)=k\frac{\partial^\alpha \nu}{\partial \vec{n}^\alpha}\quad & {\rm in}\quad\Omega,\\[2mm]
 \phantom{   (-\Delta)^\alpha  +g_n(u)}
u=0\quad & {\rm in}\quad \Omega^c
\end{array}
\end{equation}
admits a unique positive weak solution $u_{k,j,n}$ satisfying

$(i)$ the mappings $k\to u_{k,j,n}$, $j\to u_{k,j,n}$ are increasing, the mapping $n\to u_{k,j,n}$ is decreasing
\begin{equation}\label{5.6}
\mathbb{G}_\alpha[\frac{\partial^\alpha \nu}{\partial\vec{n}^\alpha}](x)-\mathbb{G}_\alpha[g_n(\mathbb{G}_\alpha
[\frac{\partial^\alpha \nu}{\partial\vec{n}^\alpha}])](x)\le u_{k,n}(x)\le
\mathbb{G}_\alpha[\frac{\partial^\alpha \nu}{\partial\vec{n}^\alpha}](x),\quad\forall x\in\Omega;
\end{equation}

$(ii)$
$u_{k,n}$ is a classical solution of (\ref{eq 2.2}).

\end{proposition}

\begin{lemma}\label{lm 5.1}
There exists $c_{58}>1$ such that
\begin{equation}\label{5.8}
 \frac1{c_{58}}\rho^\alpha(x)|x|^{-N}\le \mathbb{G}_\alpha[\frac{\partial^\alpha \delta_0}{\partial \vec{n}^\alpha}]\le c_{58}\rho^\alpha(x)|x|^{-N},\qquad
\forall x\in\Omega.
\end{equation}
\end{lemma}
{\bf Proof.} For any $x\in\Omega$, there exists $s>0$ such that $x\in \mathcal{C}_s$,
and then letting $y_t=t\vec{n}_0$ with $t\in(0,s/2)$, we have that
$$|y_t-x|=s-t> \frac s2=\frac12\max\{\rho_{\partial\Omega}(y_t),\rho_{\partial\Omega}(x)\}.$$
Thus,  it follows by apply \cite[Theorem 1.1,Theorem 1.2]{CS} that
$$\frac1{c_{59}}\frac{\rho^\alpha_{\partial\Omega}(y_t)\rho^\alpha_{\partial\Omega}(x)}{|x-y_t|^{N}}\le
G_\alpha(x,y_t) \le  c_{59}\frac{\rho^\alpha_{\partial\Omega}(y_t)\rho^\alpha_{\partial\Omega}(x)}{|x-y_t|^{N}}.$$
From
\begin{equation}
\mathbb{G}_\alpha[\frac{\partial^\alpha \delta_0}{\partial \vec{n}^\alpha}](x)=\lim_{t\to0^+}G_\alpha(x,y_t), \qquad \forall x\in\Omega,
\end{equation}
we deduce that
$$\frac1{c_{60}}\rho^\alpha(x)|x|^{-N}\le \mathbb{G}_\alpha[\frac{\partial^\alpha \delta_0}{\partial \vec{n}^\alpha}]\le c_{60}\rho^\alpha(x)|x|^{-N},\qquad
\forall  x\in\Omega.$$

\begin{lemma}\label{lemma 5.2}
Assume that   $g$ is a continuous nondecreasing function satisfying
$g(0)\ge0$, (\ref{g2}) and (\ref{g3}).
Then
 \begin{equation}\label{5.7}
 \lim_{\rho(x)\to0^+}\frac{\mathbb{G}_\alpha[g(\mathbb{G}_\alpha[\frac{\partial^\alpha \nu}{\partial \vec{n}^\alpha}])](x)}{
\rho(x)^{\alpha-1}+ \rho^\alpha(x) |x|^{-N}}=0.
 \end{equation}
\end{lemma}
{\bf Proof.} From Lemma \ref{lm 2.3} and (\ref{g2}),
 \begin{equation}\label{5.4.4}
\lim_{\rho(x)\to0^+}\mathbb{G}_\alpha[g(\mathbb{G}_\alpha[\frac{\partial^\alpha \omega}{\partial \vec{n}^\alpha}])](x)\rho^{1-\alpha}(x)=0.
\end{equation}
From \cite[Lemma 4.1]{CH}, we have that
 \begin{equation}\label{5.4.5}
\lim_{s\to0^+}\mathbb{G}_\alpha[g(\mathbb{G}_\alpha[\frac{\partial^\alpha \delta_0}{\partial \vec{n}^\alpha}])](s\vec{n}_0)s^{N-\alpha}=0.
\end{equation}
By (\ref{g3}),
$$
 \mathbb{G}_\alpha[g(\mathbb{G}_\alpha[\frac{\partial^\alpha \nu}{\partial \vec{n}^\alpha}])](x)\le \lambda\left[
 \mathbb{G}_\alpha[g(\mathbb{G}_\alpha[\frac{\partial^\alpha \omega}{\partial \vec{n}^\alpha}])](x)
+ \mathbb{G}_\alpha[g(\mathbb{G}_\alpha[\frac{\partial^\alpha \delta_0}{\partial \vec{n}^\alpha}])](x)\right],
$$
together with (\ref{5.4.4}) and (\ref{5.4.5}), we implies (\ref{5.7}).\qquad$\Box$

\medskip

\noindent{\bf Proof of Theorem \ref{teo 3}.}
 The existence of weak solution just follows the procedure of the proof of Theorem \ref{teo 1}  by using Proposition \ref{pr 4.1}
 and  Proposition \ref{pr 4.2}.
It is the same to prove the uniqueness and regularity of  weak solution.
Finally, plugging (\ref{2.1}), (\ref{5.8}) and (\ref{5.7}) into (\ref{5.6}) replaced $g_n$ by $g$, we obtain  (\ref{eq 1.3}).   \qquad$\Box$

\bigskip\medskip

\noindent{\bf Acknowledge:}  H. Chen is supported by NSFC, No: 11401270.


\begin{thebibliography}{99}
\bibitem {BM} C. Bandle, M. Marcus,
Asymptotic behaviour of solutions and derivatives for semilinear elliptic problems with blow-up on the boundary,
{\it Ann. Inst. H. Poincar¨¦ Anal. Non Lin¨¦aire, 12,} 155-171 (1995).


\bibitem {BBC} Ph. B\'{e}nilan, H. Brezis and M. Crandall, A semilinear elliptic equation in $L^1(\R^N )$, {\it Ann. Sc. Norm. Sup. Pisa Cl. Sci. 2}, 523-555 (1975).

\bibitem {CS1} L. Caffarelli and L. Silvestre, Regularity theory for fully nonlinear integro-differential equaitons,  {\it
Comm. Pure Appl. Math.}, 62(5), 597-638, 2009.

\bibitem {CS} Z. Chen and R. Song, Estimates on Green functions and poisson kernels for symmetric stable process, {\it Math.
Ann. 312}, 465-501 (1998).

\bibitem {CFQ} H. Chen, P. Felmer and A. Quaas, Large solution to elliptic  equations involving fractional Laplacian,
accepted by {\it  Ann. Inst. H. Poincar\'{e}, Analyse Non Lin\'{e}aire}, DOI: 10.1016/j.anihpc.2014.08.001.

\bibitem {CH} H. Chen and H. Hajaiej, Existence, Non-existence, Uniqueness of solutions for semilinear elliptic equations involving measures concentrated on boundary, {\it arXiv:1410.2672} (2014).

\bibitem {CC} R. Cignoli  and M. Cottlar, An Introduction to Functional Analysis, {\it North-Holland, Amsterdam}, 1974.

\bibitem {CV1} H. Chen and L. V\'{e}ron, Semilinear fractional elliptic equations involving measures, {\it J. Differential equations 257(5)}, 1457-1486 (2014).

\bibitem {DG}
Y. Du, Z. Guo,
Uniqueness and layer analysis for boundary blow-up solutions,
{\it J. Math. Pures Appl., 83,} 739-763 (2004).

\bibitem {DL} M. del Pino and R. Letelier, The influence of domain geometry in
boundary blow-up elliptic problems, {\it Nonlinear Analysis: Theory,
Methods \& Applications, 48(6),} 897-904, (2002).

\bibitem {DZZ} Y. Du, Z. Guo, F. Zhou, Boundary blow-up solutions with interior layers and spikes in a bistable problem,
{\it Discrete Contin. Dyn. Syst. 19} 271-298 (2007).

\bibitem {FQ2} P. Felmer and A. Quaas, Fundamental solutions and Liouville type theorems for nonlinear integral operators, {\it
Advances in Mathematics}, 226, 2712-2738, 2011.

\bibitem {K} J. B. Keller, On solutions of $\Delta u = f(u)$, {\it Comm. Pure Appl.
Math.},  10,  503-510, 1957.

\bibitem {GLL} J. Garcia-Meli¨¢n, R. Letelier, J. de Lis,
Uniqueness and asymptotic behaviour for solutions of semilinear problems with boundary blow-up,
{\it Proc. Amer. Math. Soc., 129}, 3593-3602 (2001)

\bibitem {GV} A. Gmira and  L. V\'{e}ron, Boundary singularities of solutions of some nonlinear elliptic equations,
{\it Duke Math. J. 64}, 271-324 (1991).

\bibitem {GZ} Z. Guo, F. Zhou,  Exact multiplicity for boundary blow-up solutions, {\it J. Differential Equations 228,} 486-506 (2006).

\bibitem {O} R. Osserman, On the inequality $\Delta u = f(u)$, {\it Pac. J. Math.} 7,
1641-1647, 1957.

\bibitem {OS}
T. Ouyang, J. Shi, Exact multiplicity of positive solutions for a class of semilinear problems,
{\it J. Differential Equations, 146,} 121-156 (1998)

\bibitem {MP} M. Pertti.  Geometry of sets and measures in Euclidean spaces: fractals and rectifiability,  {\it Cambridge University Press,} (1999).

\bibitem {MV0} M. Marcus  and L. V\'{e}ron, Uniqueness and asymptotic behavior of solutions with boundary blow-up for a class of nonlinear elliptic equations,
{\it  Ann. Inst. H. Poincar\'{e}, Analyse Non Lin\'{e}aire,   14,}   237-274 (1997).


\bibitem {MV1} M. Marcus  and L. V\'{e}ron,  Boundary trace of positive solutions to nonlinear elliptic inequalities,
{\it  Ann. Scu. Norm. Sup. Pisa}, to appear (2004).

\bibitem {MV} M. Marcus and L. V\'{e}ron, Existence and uniqueness results for large solutions of
general nonlinear elliptic equation, {\it J. Evol. Equ.} 3, 637-652 (2003).

\bibitem {TV} N. Tai  and L. V\'{e}ron, Boundary singularities of solutions to elliptic viscous Hamilton-Jacobi equations,
{\it  J. Funct. Anal. 263},  1487-1538 (2012).

\end{thebibliography}
\end{document}